\documentclass[12pt,reqno]{amsart}
\UseRawInputEncoding  

\usepackage{amsmath}
\usepackage{amsfonts}
\usepackage{amssymb}
\usepackage{caption}
\usepackage{latexsym}
\usepackage{color}
\usepackage{hyperref}
\usepackage{url}
\usepackage{enumerate}

\renewcommand{\proof}{\noindent{\it Proof.\ \ }}

 \renewcommand\a{\alpha}  \renewcommand\b{\beta}  \renewcommand\d{\delta}
 \newcommand\A{\mathrm{A}} \newcommand\AGL{\mathrm{AGL}}     \newcommand\Aut{\mathrm{Aut}}
 
\newcommand\C{\mathbf{C}}   \newcommand\Cay{\mathrm{Cay}}    
\newcommand\D{\mathrm{D}}

     \newcommand\GL{\mathrm{GL}}  
\newcommand\Ga{\mathrm{\Gamma}}

\newcommand\K{\mathsf{K}}
\newcommand\M{\mathrm{M}}  \renewcommand \mod{{\rm mod~}}

 \renewcommand\O{\mathrm{O}} 
\renewcommand\P{\mathrm{P}}  \newcommand\PGL{\mathrm{PGL}}      \newcommand\PSL{\mathrm{PSL}}     \newcommand\PSU{\mathrm{PSU}}
  
\newcommand\Q{\mathrm{Q}}

 \newcommand\SL{\mathrm{SL}}  \newcommand\soc{\mathrm{soc}}     \newcommand\Sy{\mathrm{S}}

\newcommand\ZZ{\mathbb{Z}}

\newcommand\SyG{\mathrm{S}}

\newcommand\Mult{\mathrm{Mult}}
\newcommand\diam{\mathrm{diam}}
\newcommand\BB{\mathcal{B}}
\newcommand\GG{\mathcal{G}}
\newcommand\ov{\overline}
\newcommand\la{\langle} \newcommand\ra{\rangle}

\newtheorem{theorem}{Theorem}[section]%
\newtheorem{lemma}[theorem]{Lemma}%
\newtheorem{corollary}[theorem]{Corollary}%
\newtheorem{proposition}[theorem]{Proposition}%
\newtheorem{problem}[theorem]{Problem}%
\newtheorem{example}[theorem]{Example}%
\newtheorem{question}[theorem]{Question}%

\begin{document}

\title[Two-geodesic transitive graphs of order $p^n$ with $n\leq3$]
{Two-geodesic transitive graphs of order $p^n$ with $n\leq3$}

\thanks{$^*$Corresponding author. }
\thanks{1991 MR Subject Classification 20B15, 20B30, 05C25.}
\thanks{This work was partially supported by the National Natural Science Foundation of China (11731002, 12071023, 12161141005) and the
111 Project of China (B16002).}

\author[J.-J. Huang, Y.-Q. Feng, J.-X. Zhou and F.-G. Yin]{Jun-Jie Huang, Yan-Quan Feng$^*$, Jin-Xin Zhou, Fu-Gang Yin}
\address{Jun-Jie Huang, Yan-Quan Feng, Jin-Xin Zhou, Fu-Gang Yin\\
School of mathematics and statistics\\
Beijing Jiaotong University\\
Beijing \\
100044, P. R. China}
{\email{20118006@bjtu.edu.cn(J.-J. Huang), yqfeng@bjtu.edu.cn (Y.-Q. Feng), \linebreak
jxzhou@bjtu.edu.cn (J.-X. Zhou), 18118010@bjtu.edu.cn(F.-G. Yin)}}

\begin{abstract}
A vertex triple $(u,v,w)$ of a graph is called a {\it $2$-geodesic} if $v$ is adjacent to both $u$ and $w$ and $u$ is not adjacent to $w$. A graph is said to be {\it $2$-geodesic transitive} if its automorphism group is transitive on the set of $2$-geodesics. In this paper, a complete classification of $2$-geodesic transitive graphs of order $p^n$ is given for each prime $p$ and $n\leq 3$. It turns out that all such graphs consist of three small graphs: the complete bipartite graph $\K_{4,4}$ of order $8$, the Schl\"{a}fli graph of order $27$ and its complement, and fourteen infinite families: the cycles $\C_p, \C_{p^2}$ and $\C_{p^3}$, the complete graphs $\K_p, \K_{p^2}$ and $\K_{p^3}$, the complete multipartite graphs $\K_{p[p]}$, $\K_{p[p^2]}$ and $\K_{p^2[p]}$, the Hamming graph $H(2,p)$ and its complement, the Hamming graph $H(3,p)$, and two infinite families of normal Cayley graphs on extraspecial group of order $p^3$ and exponent $p$.
\end{abstract}

\maketitle

\qquad {\textsc k}{\scriptsize \textsc {eywords.} 2-geodesic transitive graph, 2-arc transitive graph, primitive group} {\footnotesize}

\section{Introduction}
Given a graph $\Ga$, let $\Aut(\Ga)$ be the automorphism group of $\Ga$, and $\Ga$ is said to be {\em vertex-transitive}, if $\Aut(\Ga)$ is transitive on the vertex set $V(\Ga)$ of $\Ga$. If not otherwise stated, all graphs consider in this paper are finite, simple, undirected, connected and vertex-transitive. Let $u,v,w$ be three vertices of $\Ga$ such that $v$ is adjacent to both $u$ and $w$. If $u\neq w$, then the triple $(u, v, w)$ will be called a {\it $2$-arc}, and if, in addition, $u$ is not adjacent to $w$, then $(u, v, w)$ will be called a {\it $2$-geodesic}. An {\it arc} of $\Ga$ is an
ordered pair of adjacent vertices of $\Ga$. For $G\leq\Aut(\Ga)$, if $G$ is transitive on the set of arcs of $\Ga$, then $\Ga$ is said to be {\it $G$-arc transitive}, and if $G$ is transitive on the set of $2$-arcs of $\Ga$, then $\Ga$ is said to be {\it $(G,2)$-arc transitive}. Similarly, if $G$ is transitive on the set of $2$-geodesics of $\Ga$, then $\Ga$ is said to be {\it $(G,2)$-geodesic transitive}. When $G=\Aut(\Ga)$, a $G$-arc transitive graph,  $(G,2)$-arc transitive graph or $(G,2)$-geodesic transitive graph will be simply called an {\it arc transitive graph}, a {\it $2$-arc transitive graph} or {\it $2$-geodesic transitive graph}, respectively. Clearly, $2$-arc transitive graphs are also $2$-geodesic transitive. Following the seminal work of Tutte~\cite{Tutte1947}, the $2$-arc transitive graphs have been received considerable attention, see \cite{ACX,DMO,Li,LSS,MP,Praeger92} and references therein. Note that if $\Ga$ has girth $3$, then the $2$-arcs contained in $3$-cycles are not $2$-geodesics. There exist $2$-geodesic transitive but not $2$-arc transitive graphs, and a typical example is the octahedron. The $2$-geodesic transitive but not $2$-arc transitive graphs have also been extensively studied in the literature, see, for example, \cite{DJLP13,DJLP14,DJLP15,Jin17-2,Jin19,JDLP,JLW}.

There are a fair amount of works on $2$-arc transitive graphs with particular orders in the literature. One of the remarkable results is the classification of basic $2$-arc transitive graphs of prime-power order (\cite{Li}). Later, this was extended to the class of vertex transitive locally primitive graphs of prime-power order in \cite{LMP2009}. All the $2$-arc transitive graphs of order a prime or a prime square have also been classified in \cite{MP}. Recently, Jin~\cite{Jin17} initiated the study of $2$-geodesic transitive graphs of prime power order, and he proved that for a $2$-geodesic transitive graph $\Ga$ of prime power, if $\Aut(\Ga)$ is quasiprimitive on $V(\Ga)$, then $\Ga$ is isomorphic to some well-known graphs (see Proposition~\ref{quas-2-geo}), or $\Aut(\Ga)$ is primitive on $V(\Ga)$ of affine type. The following problem was proposed in \cite{Jin17}.\medskip

\begin{problem}{\rm(\cite[Problem 1.4]{Jin17})}\label{problem}
Let $\Ga$ be a $2$-geodesic transitive graph of prime power order which
is not $2$-arc transitive. Classify such graphs where $\Aut(\Ga)$ acts quasiprimitively
on $V(\Ga)$ of affine type.
\end{problem}

Inspired by the facts listed above, in this paper we shall focus on the classification of 2-geodesic transitive graphs of order $p^n$ with $p$ a prime and $n\leq 3$. The following is our first main result which provides a partial solution to Problem~\ref{problem}.

\begin{theorem}\label{prim-p23}
Let $p$ be a prime and let $\Ga$ be a connected $2$-geodesic transitive but not $2$-arc transitive normal Cayley graph on $\ZZ_p^n$ with $n\leq 3$. Then $\Aut(\Ga)$ is primitive on $V(\Ga)$ if and only if $\Ga\cong H(2,3)$ or $H(3,3)$.
\end{theorem}

Another motivation for us to consider 2-geodesic transitive graphs of order $p^n$ with $p$ a prime and $n\leq 3$ comes from the work in \cite{CJL,HFZ}. Before proceeding, we introduce the following definitions. For two vertices $u$ and $v$ of a graph $\Ga$, the {\it distance} $d_\Ga(u,v)$ between $u$ and $v$ in $\Ga$ is the smallest length of paths between $u$ and $v$, and the {\it diameter} $\diam(\Ga)$ of $\Ga$ is the maximum distance occurring over all pairs of vertices. A graph $\Ga$ is called {\it distance transitive} if, for any vertices $u,v,x,y$ with $d_\Ga(u,v)=d_\Ga(x,y)$, $\Ga$ has an automorphism mapping the pair $(u,v)$ to $(x,y)$. For a positive integer $s$, $\Ga$ is said to be {\it $s$-distance transitive} if, for each $1\leq i\leq s$, $\Aut(\Ga)$ is transitive on the ordered pairs of form $(u,v)$ with $d_\Ga(u,v)=i$. In \cite{CJL}, Chen et al. classified $2$-distance transitive circulants, and they proposed the following question:

\begin{question}{\rm(\cite[Question 1.2]{CJL})}\label{question}
Is there a normal Cayley graph which is $2$-distance transitive,
but neither distance transitive nor $2$-arc transitive?
\end{question}

In \cite{HFZ}, we answered this question in positive by constructing an infinite family of 2-distance transitive Cayley graphs of order a prime cube which are neither distance transitive nor 2-arc transitive. This family of graphs are also $2$-geodesic transitive. Stimulated by this, in this paper, we give a complete classification of the $2$-geodesic transitive graphs of order $p^n$ with $p$ a prime and $n\leq 3$.  As a result, we find a new family of $2$-geodesic transitive Cayley graphs which are not $2$-arc transitive.

Note that the $2$-geodesic transitive Cayley graphs on cyclic groups were classified in \cite{JLW}, from which we find that a $2$-geodesic transitive graph of order a prime is either a cycle or a complete graph. In our next two main theorems, we give a complete classification of $2$-geodesic transitive graphs of order $p^2$ or $p^3$ with $p$ a prime.

\begin{theorem}\label{Thm-p^2} Let $p$ be a prime and let $\Ga$ be a connected arc-transitive graph of order $p^2$. Then $\Ga$ is $2$-geodesic transitive if and only if one of the following holds:
\begin{itemize}
  \item [(1)] If $\Ga$ is $2$-arc transitive, then $\Ga$ is isomorphic to $\C_{p^2}$ or $\K_{p^2}$;
  \item [(2)] If $\Ga$ is not $2$-arc transitive, then $\Ga$ is isomorphic to $\K_{p[p]}$, the Hamming graph $H(2,p)$ or its complement $\ov{H(2,p)}$, where $p\geq 3$.
\end{itemize}
\end{theorem}

\begin{theorem}\label{Thm-p^3}
Let $p$ be a prime and let $\Ga$ be a connected arc transitive graph of order $p^3$. Then $\Ga$ is $2$-geodesic transitive if and only if one of the following holds:
\begin{itemize}
  \item [(1)] If $\Ga$ is $2$-arc transitive, then $\Ga$ is isomorphic to $H(3,2)$, $\K_{4,4}$, $\C_{p^3}$ or $\K_{p^3}$;
  \item [(2)] If $\Ga$ is not $2$-arc transitive, then $\Ga$ is isomorphic to one of the following graphs:
  \begin{itemize}
    \item[(i)] the Schl\"{a}fli graph or its complement;
    \item[(ii)] $\K_{p^2[p]}$ with $p\geq2$, or $\K_{p[p^2]}$ with $p\geq3$;
    \item[(iii)] the Hamming graph $H(3,p)$ with $p\geq3$;
    \item[(iv)] the normal Cayley graph $\Cay(G,S)$, where $G=\la a,b,c\mid a^p=b^p=c^p=1,[a,b]=c,[a,c]=[b,c]=1\ra$ with $p\geq3$, and $S=\{a^i,b^i\mid i\in\ZZ_p^*\}$ or $S=\{a^i,b^i, (b^jab^j)^i\mid i,j\in\ZZ_p^*\}$, respectively.
  \end{itemize}
\end{itemize}
\end{theorem}

\noindent{\bf Remark:}\ All graphs in Theorems~\ref{prim-p23}, \ref{Thm-p^2} and \ref{Thm-p^3} will be explicitly defined in Sections~\ref{Example} and \ref{Newfamily}, and all but graphs in part (iv) of Theorem~\ref{Thm-p^3}~(2) have been well-studied in the literature. In part (iv) of Theorem~\ref{Thm-p^3}~(2), the graphs $\Cay(G,S)$ with $S=\{a^i,b^i\mid i\in\ZZ_p^*\}$ were first given in \cite{HFZ} while the graphs $\Cay(G,S)$ with $S=\{a^i,b^i, (b^jab^j)^i\mid i,j\in\ZZ_p^*\}$ are a new family of $2$-geodesic transitive but not 2-arc transitive graphs.

\medskip

The layout of this paper is as follows. In Section~\ref{Preliminaries}, we list some definitions and some preliminary results, and in Sections~\ref{Example} and \ref{Newfamily}, all $2$-geodesic transitive graphs in question are listed, where a new infinite family of $2$-geodesic transitive Cayley graphs on the extraspecial group of order $p^3$ and exponent $p$ is constructed in Section~\ref{Newfamily}. In Section~\ref{Affine case}, we prove Theorem~\ref{prim-p23}, and in Section~\ref{MainProof}, we prove Theorems~\ref{Thm-p^2} and \ref{Thm-p^3}.

\section{Preliminaries}\label{Preliminaries}

In this section, we give some results which will be used in following sections. Let $p$ be a prime and let $n$ be a positive integer. We denote by $n_p$ the largest $p$-power dividing $n$, and by $n_{p'}$ the largest divisor of $n$ such that $p\nmid n_{p'}$. We use $\ZZ_n$ and $\ZZ_p^r$ to denote the cyclic group of order $n$ and the elementary abelian group of order $p^r$, respectively. Also, we denote by $\ZZ_p$ the additive group of integers modulo $p$, and by $\ZZ_p^*$ the multiplicative group of all non-zero numbers in $\ZZ_p$. Denote by $\Q_8$, $\D_{2n}$, $\A_n$ and $\SyG_n$, the quaternion group of order $8$, the dihedral group of order $2n$, the alternating group of degree $n$, and the symmetric group of degree $n$, respectively. For a subgroup $H$ of a group $G$, denote by $H^*$ the subset of $H$ deleting $1$, that is, $H^*=\{h\ |\ h\in H, h\not=1\}$. For two groups $A$ and $B$, denote by $A\times B$ the direct product of $A$ and $B$, by $A\rtimes B$ or $A:B$ a semidirect product of $A$ by $B$, and by $A.B$ an extension of $A$ by $B$. For a group $G$, denote by $Z(G)$ the central of $G$, and by $\soc(G)$ the socle of $G$, that is, the product of all minimal normal subgroup of $G$.

For graph notation, denote by $\C_n$ the cycle of length $n$, $\K_n$ the complete graph of order $n$, $\K_{n,n}$ the complete bipartite graph of order $2n$, $\K_{n,n}-n\K_2$ the subgraph of $\K_{n,n}$ minus a matching, and $\K_{m[b]}$ the complete multipartite graph consisting $m\geq3$ parts of size $b\geq2$, respectively. Let $\Ga$ be a graph.
Denote by $\ov{\Ga}$ the complement graph of $\Ga$, that is, the graph with the vertex set $V(\Ga)$ and two vertices are adjacent in $\ov{\Ga}$ if and only if they are not adjacent in $\Ga$. A {\it clique} of a graph $\Ga$ is a complete subgraph and a {\it maximal clique} is a clique which is not contained in a larger clique. For two disjoint subsets $U,W\subset V(\Ga)$, denote by $[U]$ the subgraph of $\Ga$ induced by $U$, and by $[U,W]$ the subgraph of $\Ga$ with vertex set $U\cup W$ and edge set consisting of the edges of $\Ga$ between $U$ and $W$.
Moreover, we use $\Ga_i(u)$, with $i\leq\diam(\Ga)$ and $u\in V(\Ga)$, to denote the set of vertices with distance $i$ from $u$ in $\Ga$.

A classification of finite groups of order prime cube was given in \cite[P. 8, Theorem 2.14]{Gore} and \cite[P. 203, Theorem 5.1]{Gore}.

\begin{proposition}\label{group-p^3}
 Let $p$ be an odd prime and let $G$ be a finite group of order $p^3$. Then $G$ is isomorphic to one of the following groups: $\ZZ_{p^3},\ZZ_{p^2}\times\ZZ_p,\ZZ_p^3,\la a,b\mid a^{p^2}=b^p=1,a^b=a^{1+p}\ra$ or $\la a,b,c\mid a^p=b^p=c^p=1,[a,b]=c,[a,c]=[b,c]=1\ra$.
\end{proposition}

All non-abelian simple groups with a subgroup of prime-power index were determined in \cite[Theorem 1]{Gur}, and we record it here.

\begin{proposition}\label{simple-ind}
Let $T$ be a nonabelian simple group with a subgroup $H$ of index $p^n$, where $p$ is a prime and $n\geq1$. Then one of the following holds:
\begin{itemize}
  \item[(1)] $T\cong\A_{p^n}$ and $H\cong\A_{p^n-1}$;
  \item[(2)] $T\cong\PSL(d,q)$ and $H$ is a maximal parabolic subgroup of $T$, with $p^n=(q^d-1)/(q-1)$ and $d$ a prime;
  \item[(3)] $T\cong\PSL(2,11)$ and $H\cong\A_5$, with $p^n=11$;
  \item[(4)] $T\cong\M_{11}$ and $H\cong\M_{10}$, with $p^n=11$;
  \item[(5)] $T\cong\M_{23}$ and $H\cong\M_{22}$, with $p^n=23$;
  \item[(6)] $T\cong\PSU(4,2)$ and $H\cong\ZZ_2^4\rtimes\A_5$, with $p^n=27$.
\end{itemize}
\end{proposition}

For a group $G$,  $M\lessdot G$ means that $M$ is a maximal subgroup of $G$. All maximal subgroups of $\PSL(2,q)$ and $\PGL(2,q)$ were given in \cite[P. 392-418]{Suzu} and \cite[Theorem 2]{COT}.

\begin{proposition}\label{maxPGL2}
Let $p\geq5$ be a prime.
\begin{itemize}
  \item[(1)] If $M\lessdot \PGL(2,p)$, then $M$ is one of the following groups: $\PSL(2,p)$, $\ZZ_p\rtimes\ZZ_{p-1}$, $\D_{2(p-1)}$ with $p\neq 5$, $\D_{2(p+1)}$, or $\SyG_4$ with $p\equiv\pm3(\mod 8)$;
  \item[(2)] If $M \lessdot \PSL(2,p)$, then $M$ is one of the following groups: $\ZZ_p\rtimes\ZZ_{(p-1)/2}$, $\D_{p-1}$ with $p\geq 13$, $\D_{p+1}$ with $p\neq 7$, $\A_4$ with $p\equiv\pm3(\mod 8)$ and
      $p\equiv\pm3(\mod 10)$,
      $\SyG_4$ with $p\equiv\pm1(\mod 8)$, or $\A_5$ with $p\equiv\pm1(\mod 10)$.
\end{itemize}
\end{proposition}

By Proposition~\ref{maxPGL2}, the smallest index of a maximal subgroup of $\PSL(2,p)$ is $p$, and so every non-trivial transitive action of $\PSL(2,p)$ has degree at least $p$.

The following can be read off from \cite[Tables 8.3 and 8.4]{BHR}, where the maximal subgroups of $\PGL(3,p)$ corresponding to row~1,  row~4 and row~5 of \cite[Table 8.3]{BHR} should be  $\ZZ_p^2:\GL(2,p)$, $(\ZZ_{p-1}\times\ZZ_{p-1}):\SyG_3$ and $\ZZ_{p^2+p+1}:\ZZ_3$
(see \cite[Chapter 4]{KL}), respectively.
The notation $\frac{1}{d}\GL(2,p)$ is the factor group of $\GL(2,p)$ by the subgroup of order $d$ in $Z(\GL(2,p))$.

\begin{proposition}\label{maxPGL3}
Let $p\geq5$ be a prime and let $d=(3,p-1)$.
\begin{itemize}
  \item[(1)] If $M\lessdot \PGL(3,p)$, then $M$ is one of the following groups: $\PSL(3,p)$, $\ZZ_p^2:\GL(2,p)$, $(\ZZ_{p-1}\times\ZZ_{p-1}):\SyG_3$, $\ZZ_{p^2+p+1}:\ZZ_3$, and $(\ZZ_3^2:\Q_8).\ZZ_3$ with $p\equiv1(\mod 3)$ and $p\not\equiv1(\mod 9)$;

   \item[(2)]If $M\lessdot \PSL(3,p)$, then $M$ is one of the following groups: $\ZZ_p^2:\frac{1}{d}\GL(2,p)$, $(\ZZ_{(p-1)/d}\times \ZZ_{p-1}):\SyG_3$, $\ZZ_{(p^2+p+1)/d}:\ZZ_3$, $(\ZZ_3^2:\Q_8).\ZZ_{\frac{(p-1,9)}{3}}$ with $p\equiv1(\mod 3)$, $\PGL(2,p)$, $\PSL(2,7)$ with $p\equiv1,2,4(\mod 7)$, and $\A_6$ with $p\equiv1,4 (\mod 15)$.

\end{itemize}
\end{proposition}

By Proposition~\ref{maxPGL3}, the smallest index of a maximal subgroup of $\PSL(3,p)$ is $p^2+p+1$, and so every non-trivial transitive action of $\PSL(3,p)$ has degree at least $p^2+p+1$.

From Propositions~\ref{maxPGL2} and \ref{maxPGL3}, all $2$-transitive actions of subgroups of $\PGL(3,p)$ can be determined.

\begin{lemma}\label{2-transPGL3}
Let $p\geq 5$ be a prime.
Assume that $X\leq \PGL(3,p)$ has a $2$-transitive action on $\Omega$ with $n=|\Omega|\geq 4$. Let $K$ be the kernel of $X$ on $\Omega$.
If $\soc(X/K)$ is elementary abelian, then either
\begin{itemize}
    \item[(1)] $(n,X/K)=(4,\A_4)$, $(4,\SyG_4)$, $(p,\AGL(1,p))$, or $(9,L)$ with $L\leq (\ZZ_3^2:\Q_8).\ZZ_3$ and $1\equiv(\mod 3)$; or
    \item[(2)] $n=p^2$ and $X/K\leq \AGL(2,p)$.
\end{itemize}

If $\soc(X/K)$ is nonabelian simple, then one of the following holds:
\begin{itemize}
  \item[(1)] $(n,X/K)=(5,\A_5)$, $(6,\A_6)$, $(10,\A_6)$, $(8,\PSL(2,7))$ or $(11,\PSL(2,11))$, and if $X/K=\A_6$,  $\PSL(2,7)$ or $\PSL(2,11)$
       then $p\equiv 1,4(\mod 15)$, $p\equiv 1,2,4 (\mod 7)$ or $p=11$, respectively;
    \item[(2)] $n=p+1$, and $X/K=\PSL(2,p)$ or $\PGL(2,p)$;
    \item[(3)] $n=p^2+p+1$, and $X/K=\PSL(3,p)$ or $\PGL(3,p)$.
    \end{itemize}
\end{lemma}

\proof By \cite[Theorem 4.1B]{Dixon}, either $\soc(X/K)\cong T$, a non-abelian simple group, or $\soc(X/K)\cong \ZZ_q^r$ with $q$ a prime and $r$ a positive integer. For $\soc(X/K)\cong T$,
by Propositions~\ref{maxPGL2} and \ref{maxPGL3} we have $T=\A_5$, $\A_6$, $\PSL(2,p)$ or $\PSL(3,p)$, and then all $2$-transitive actions of $T$ follow from \cite[Table 2.1]{LSS}. Furthermore, $X/K=\A_5$ or $\A_6$ for $T=\A_5$ or $\A_6$, $X/K=\PSL(2,p)$ or $\PGL(2,p)$ for $T=\PSL(2,p)$, and $X/K=\PSL(3,p)$ or $\PGL(3,p)$ for $T=\PSL(3,p)$.

Let $\soc(X/K)\cong\ZZ_q^r$. Then $n=|\Omega|=q^r$. If $n=4$, then $X/K=\A_4$ or $\SyG_4$. We further assume $n\geq 5$. First we prove a claim.

\medskip
\noindent{\bf Claim:} Assume $Z_1\leq Z(\GL(2,p))$ and $\GL(2,p)/Z_1$ has a subgroup isomorphic to $X/K_1$ for some $K_1\leq K$. Then $n=p$ and $X/K\cong\AGL(1,p)$.

Set $\ov{\GL(2,p)}=\GL(2,p)/Z_1$ and $\ov{Z(\GL(2,p))}=Z(\GL(2,p))/Z_1$. Then $\ov{Z(\GL(2,p))}\leq Z(\ov{\GL(2,p)})$ and $\ov{\GL(2,p)}/\ov{Z(\GL(2,p))}\cong\PGL(2,p)$.
If $\ov{Z(\GL(2,p))}\not=Z(\ov{\GL(2,p)})$, then $Z(\ov{\GL(2,p)})/\ov{Z(\GL(2,p))}$ is  a non-trivial solvable normal subgroup of $\ov{\GL(2,p)}/\ov{Z(\GL(2,p))}$, which is impossible because $\soc(\PGL(2,p))=\PSL(2,p)$. It follows that $\ov{Z(\GL(2,p))}=Z(\ov{\GL(2,p)})$ and hence $\ov{\GL(2,p)}/Z(\ov{\GL(2,p)})\cong\PGL(2,p)$.

Set $\ov{X}=X/K_1$ and $\ov{K}=K/K_1$. Then $X/K\cong \ov{X}/\ov{K}$ and $\soc(\ov{X}/\ov{K})\cong \soc(X/K)\cong \ZZ_q^r$. Since $K_1\leq K$, $\ov{X}$ has a $2$-transitive action on $\Omega$ with kernel $\ov{K}$, and hence $\ov{X}/\ov{K}$ is a $2$-transitive permutation group on $\Omega$ with the unique regular subgroup $\soc(\ov{X}/\ov{K})\cong \ZZ_q^r$.

Suppose $Z(\ov{X}/\ov{K})\not=1$. The $2$-transitivity of $\ov{X}/\ov{K}$ implies that every non-trivial subgroup of $Z(\ov{X}/\ov{K})$ is transitive on $\Omega$, forcing $|\Omega|=q$ and $Z(\ov{X}/\ov{K})\cong \ZZ_q$. It follows that $\ov{X}/\ov{K}\cong \ZZ_q\rtimes \ZZ_{q-1}$, and since  $Z(\ov{X}/\ov{K})\cong \ZZ_q$, we have $\ov{X}/\ov{K}\cong \ZZ_q\times \ZZ_{q-1}$, contradicting the $2$-transitivity of $\ov{X}/\ov{K}$.

Thus, $Z(\ov{X}/\ov{K})=1$, and so $Z(\ov{X})\leq \ov{K}$. By hypothesis, without loss of any generality, we assume $X/K_1\leq \GL(2,p)/Z_1$. Then $\ov{X}=X/K_1\leq \GL(2,p)/Z_1=\ov{\GL(2,p)}$. Furthermore,
$\ov{X}\cap Z(\ov{\GL(2,p)})\leq Z(\ov{X})\leq \ov{K}$,
and $\ov{X}/(\ov{X}\cap Z(\ov{\GL(2,p)}))$ has a $2$-transitive action on $\Omega$. Since $n=|\Omega|\geq 5$, $\ov{X}/(\ov{X}\cap Z(\ov{\GL(2,p)}))$ cannot be abelain or dihedral, and since
$$\ov{X}/(\ov{X}\cap Z(\ov{\GL(2,p)}))\cong \ov{X}Z(\ov{\GL(2,p)})/Z(\ov{\GL(2,p)})\leq \ov{\GL(2,p)}/Z(\ov{\GL(2,p)})\cong \PSL(2,p),$$
by Proposition~\ref{maxPGL2} we obtain $\ov{X}/(\ov{X}\cap Z(\ov{\GL(2,p)}))\cong\ZZ_p:\ZZ_{p-1}$, which implies that  $\ov{X}/(\ov{X}\cap Z(\ov{\GL(2,p)}))$ is faithful on $\Omega$, that is, $\ov{X}\cap Z(\ov{\GL(2,p)})=\ov{K}$. Thus, $X/K\cong \ov{X}/\ov{K}\cong \ZZ_p:\ZZ_{p-1}$, and the $2$-transitivity of $X/K$ implies $X/K\cong \AGL(1,p)$, as claimed.

\medskip

Since $X\leq \PGL(3,p)$, we have $X\leq M$ with $M\lessdot \PGL(3,p)$ or $M\lessdot \PSL(3,p)$. We argue by contradiction and we suppose that $(n,X/K)\not=(9,\ZZ_3^2:\Q_8)$, $(9,(\ZZ_3^2:\Q_8).\ZZ_3)$, $(p,\AGL(1,p)$, or $(p^2, L)$ with $L\leq \AGL(2,p)$.
By Proposition~\ref{maxPGL3}, $M=\PGL(2,p)$, $\ZZ_p^2:\GL(2,p)$, $\ZZ_p^2:\frac{1}{d}\GL(2,p)$, $\A_6$, $(\ZZ_{p-1}\times\ZZ_{p-1}):\SyG_3$, $(\ZZ_{(p-1)/d}\times \ZZ_{p-1}):\SyG_3$, $\ZZ_{p^2+p+1}:\ZZ_3$ or $\ZZ_{(p^2+p+1)/d}:\ZZ_3$. We deal with these case by case.

Assume $X\leq \PGL(2,p)$. By taking $K_1=1$ and $Z_1=Z(\GL(2,p))$, Claim implies that $(n,X/K)=(p,\AGL(1,p))$, a contradiction.

Assume $M=\ZZ_p^2:\GL(2,p)$ or $\ZZ_p^2:\frac{1}{d}\GL(2,p)$. If $X\cap \ZZ_p^2\not\leq K$, then $K(X\cap \ZZ_p^2)/K$ is the unique regular subgroup of $X/K$ on $\Omega$, and since $K(X\cap \ZZ_p^2)/K\cong \ZZ_p$ or $\ZZ_p^2$, either $(n,X/K)=(p,\AGL(1,p))$ or $(n,X/K)=(p^2,L)$ with $L\leq \AGL(2,p)$, a contradiction. Thus, $X\cap \ZZ_p^2\leq K$, and $X/(X\cap \ZZ_p^2)\cong X\ZZ_p^2/\ZZ_p^2\leq \GL(2,p)$ or  $\frac{1}{d}\GL(2,p)$.
By taking $K_1=X\cap \ZZ_p^2$ and $Z_1=1$ or the subgroup of order $d$ of $Z(\GL(2,p))$ for $\GL(2,p)$ or $\frac{1}{d}\GL(2,p)$ respectively, Claim implies that $(n,X/K)=(p,\AGL(1,p))$, a contradiction.

Assume $M=\A_6$. By Atlas~\cite{Atlas}, $X\leq \A_4$, which is impossible because $n\geq 5$.

Assume $M=\ZZ_{p^2+p+1}:\ZZ_3$. Let $N=\ZZ_{p^2+p+1}$. Since $X$ is not abelian, $X\cap N\not=X$, and since $X/(X\cap N)\cong XN/N\leq M/N\cong \ZZ_3$, we have $X=(X\cap N).\ZZ_3$. Thus, $(X\cap N)K=X$ or
$(X\cap N)K=X\cap N$. For the former, $X/K$ is abelian, a contradiction. For the latter, $K\leq X\cap N$, and $(X\cap N)/K$ is the regular subgroup of the permutation group $X/K$ on $\Omega$. Since  $(X/K)/((X\cap N)/K)\cong X/(X\cap N)\leq \ZZ_3$, we have $X/K=(X\cap N).\ZZ_3$, forcing $n=4$, a contradiction.
A similar argument gives rise to a contradiction for $M=\ZZ_{(p^2+p+1)/d}:\ZZ_3$.

Assume $M=(\ZZ_{p-1}\times \ZZ_{p-1}):\SyG_3$. Write $Y=X\cap (\ZZ_{p-1}\times \ZZ_{p-1})$. Since $Y$ is abelian, $Y\not=X$, and since $X/Y\cong X(\ZZ_{p-1}\times \ZZ_{p-1})/(\ZZ_{p-1}\times \ZZ_{p-1})\leq M/(\ZZ_{p-1}\times \ZZ_{p-1})\cong \SyG_3$, we have
$Y<X\leq Y.\SyG_3$.
Since $n\geq 5$, $Y\not\leq K$, and so $KY/K$ is the regular subgroup of $X/K$. Since $(X/K)/(KY/K)\cong X/KY$ and
$X/Y\leq\SyG_3$, we have
$X/K=(KY/K).\ZZ_2$, $(KY/K).\ZZ_3$ or $(KY/K).\SyG_3$, of witch the first two are impossible. Thus, $X/K=(KY/K).\SyG_3$, and since $\soc(X/K)=KY/K\cong\ZZ_q^r$, the $2$-transitivity of $X/K$ on $\Omega$ implies that $(q^r-1)\mid 6$. Since $n=q^r\geq 5$, we have $q^r=7$ and $X/K\cong \ZZ_7\rtimes\ZZ_6$,
which implies that $(X/K)/(KY/K)\cong \ZZ_6$, a contradiction. A similar argument gives rise to a contradiction for $M=(\ZZ_{(p-1)/d}\times \ZZ_{p-1}):\SyG_3$. This completes the proof. \qed

An extension $G=N.H$ is called a {\it central extension} if $N$ is a central subgroup of $G$, and a central extension $G=N.H$ is called a {\it covering group} of $H$ if $G$ is perfect, that is, the derived group $G'$ of $G$ is equal to $G$. Schur~\cite{Schur} showed that a simple group $T$ possesses a universal covering group $G$ with the property that every covering group of $T$ is a homomorphic image of $G$, and the center $Z(G)$ is called the {\it Schur multiplier} of $T$, denoted by $\Mult(T)$. The following result was given in \cite[Lemma 2.11]{PLHL}.

\begin{proposition}\label{mult}
Let $N$ be a group of order prime or prime-square, and let $T$ be a nonabelian simple group. Then every extension $G=N.T$ is a central extension. Furthermore, $G=NG'$ and $G'= M.T$ with $M\leq N\cap\Mult(T)$.
\end{proposition}

For a finite group $G$ and an inverse closed subset $S\subseteq G \setminus\{1\}$, the {\it Cayley graph} $\Cay(G,S)$ on $G$ with respect to $S$ is defined to be the graph with vertex set $G$ and edge set $\{\{g,sg\}\mid g\in G,s\in S\}$. For $g\in G$, let $R(g)$ be the permutation of $G$ defined by $x\mapsto xg$ for all $x\in G$, and then $R(G):=\{R(g)\mid g\in G\}$ is a regular subgroup of automorphisms of  $\Cay(G,S)$. A Cayley graph $\Cay(G,S)$ is called {\em normal} if $R(G)$ is a normal subgroup of $\Aut(\Cay(G,S))$. It is widely known that a graph $\Ga$ is a Cayley graph on $G$ if and only if $\Aut(\Ga)$ has a regular subgroup isomorphic to $G$, see \cite{Sabidussi}. Moreover, $\Aut(G,S):=\{\a\in\Aut(G)\mid S^\a=S\}$ is also a subgroup of $\Aut(\Ga)$, which fixes the vertex $1$. Godsil~\cite{Godsil} proved that $R(G)\rtimes\Aut(G,S)$ is the normalizer of $R(G)$ in $\Aut(\Cay(G,S))$, and together with Xu~\cite{Xu98}, we have the following proposition.

\begin{proposition}\label{Aut-Cay}
Let $\Ga=\Cay(G,S)$ be a Cayley graph on a finite group $G$ with respect to $S$, and let $A=\Aut(\Ga)$. Then $N_A(R(G))=R(G)\rtimes\Aut(G,S)$, and $\Ga$ is a normal Cayley graph if and only if $A_1=\Aut(G,S)$, where $A_1$ is the stabilizer of $1$ in $A$.
\end{proposition}

A Cayley graph $\Cay(G,S)$ is called a {\it circulant} if $G$ is a cyclic group. The 2-arc transitive and 2-geodesic transitive circulants were classified in \cite{ACX,JLW}, respectively. We summarize those results as follows:

\begin{proposition}\label{2-geo-cir}
Let $\Ga$ be a connected $2$-geodesic transitive circulant. Then
\begin{itemize}
  \item[(1)] If $\Ga$ is $2$-arc transitive, then $\Ga$ is one of the following graphs: $\K_n$ with $n\geq 1$, $\C_n$ with $n\geq 4$, $\K_{\frac{n}{2},\frac{n}{2}}$ with $n\geq 6$, and $\K_{\frac{n}{2},\frac{n}{2}}-\frac{n}{2}\K_2$ with $\frac{n}{2}\geq5$ odd;
  \item[(2)] If $\Ga$ is not $2$-arc transitive, then $\Ga\cong\K_{m[b]}$ for some $m\geq3$ and $b\geq2$.
\end{itemize}
\end{proposition}

A Cayley graph $\Cay(G,S)$ is called a {\it normal $(X,2)$-geodesic transitive} if it is $(X,2)$-geodesic transitive for a group $X$ such that $R(G)\leq X\leq R(G)\rtimes\Aut(G,S)$.

\begin{proposition}\label{nor-2-geo}{\rm(\cite[Theorem 1.2]{DJLP14})}
Let $\Ga=\Cay(G,S)$ be a connected normal $(X,2)$-geodesic transitive Cayley graph. Then one of the following holds:
\begin{itemize}
  \item[(1)] $\Ga\cong\C_r$ and $G\cong\ZZ_r$ for some $r\geq4$;
  \item[(2)] $\Ga\cong\K_{4[2]}$ and $G\cong\Q_8$, the quaternion group, with $S=G\setminus Z(G)$;
  \item[(3)] There is a prime $q$ and an integer $m$ such that for all $a\in S$, $a$ has order $q$ with $\la a\ra^*\subseteq S$ and $\la a\ra^*\not=S$, and $b$ has order $m$ for each $b\in S^2\setminus(S\cup\{1\})$.
\end{itemize}
\end{proposition}

Let $\Ga$ be a graph with $G\leq \Aut(\Ga)$, and let $N$ be a normal subgroup of $G$ such that $N$ is intransitive on $V(\Ga)$. The {\it normal quotient graph} $\Ga_N$ of $\Ga$ induced by $N$ is defined to be the graph with vertex set the orbits of $N$ and with two orbits adjacent in $\Ga_N$ if $\Ga$ has an edge incident to some vertices in the two orbits, respectively. Furthermore, $\Ga$ is called a {\it normal $N$-cover} or simply {\it cover} of $\Ga_N$ if for every vertex $v\in V(\Ga)$ and the orbit $O(v)$ of $N$ containing $v$, $v$ and $O(v)$ has the same valency in $\Ga$ and $\Ga_N$, respectively.
The graph $\Ga$ is called {\it locally $(G,s)$-distance transitive} with $s\geq 1$, if for every vertex $u\in V(\Ga)$, $G_u$ is transitive on $\Ga_i(u)$ for all $i\leq s$, where $\Ga_i(u)$ is the set of vertices with distance $i$ from $u$ in $\Ga$. The following proposition gives a reduction for studying locally $(G,s)$-distance transitive graphs.

\begin{proposition}\label{redu}{\rm(\cite[Lemma 5.3]{DGLP12})}
Let $\Ga$ be a connected locally $(G,s)$-distance transitive graph with $s\geq2$. Let $1\neq N\lhd G$ be intransitive on $V(\Ga)$, and let $\BB$ be the set of $N$-orbits on $V(\Ga)$. Then one of the following holds:
\begin{itemize}
  \item [(1)] $|\BB|=2$;
  \item [(2)] $\Ga$ is bipartite, $\Ga_N\cong\K_{1,r}$ with $r\geq2$ and $G$ is intransitive on $V(\Ga)$;
  \item [(3)] $s=2$, $\Ga\cong\K_{m[b]},\Ga_N\cong\K_m$, where $m\geq3$ and $b\geq 2$;
  \item [(4)] $N$ is semiregular on $V(\Ga)$, $\Ga$ is a cover of $\Ga_N$, $|V(\Ga_N)|<|V(\Ga)|$ and $\Ga_N$ is locally $(G/N,s')$-distance transitive, where $s'=\min\{s,\diam(\Ga_N)\}$.
\end{itemize}
\end{proposition}

The next result illustrates that a connected locally $(G,s)$-distance transitive with $s\geq2$, cannot be a cover of  $\K_{m[b]}$.

\begin{proposition}\label{cover}{\rm(\cite[Proposition 4.2]{DGLP12})}
Let $\Ga$ be a connected locally $(G,s)$-distance transitive graph with $s\geq2$. Then there exists no nontrivial $N\lhd G$ such that $\Ga$ is a cover of $\Ga_N$ and $\Ga_N\cong\K_{m[b]}$ for some $m\geq3$ and $b\geq 2$.
\end{proposition}

A transitive permutation group $G$ on $\Omega$ is called {\em primitive} if it has only trivial blocks in $\Omega$ (a
{\em block}
is a non-empty subset $\Delta$ of $\Omega$ such that $\Delta^g\cap \Delta\not=\emptyset$ implies $\Delta^g=\Delta$ for every $g\in G$ and a block $\Delta$ is
{\em trivial}
if $|\Delta|=1$ or $|\Delta|=|\Omega|$), and {\em quasiprimitive} if every non-trivial normal subgroup of $G$ is transitive on $\Omega$. Clearly, primitive group is  quasiprimitive, but the converse is not true. Combining \cite[Theorem 1.2]{Jin17} and \cite[Theorem 1.3]{Jin21}, we have the following result.

\begin{proposition}\label{quas-2-geo}
For a prime $p$ and a positive integer $n$, let $\Ga$ be a connected $2$-geodesic but not $2$-arc transitive graph of order $p^n$, and let $\Aut(\Ga)$ be quasiprimitive on $V(\Ga)$. Then  $\Aut(\Ga)$ is primitive on $V(\Ga)$ and one of the following holds:
\begin{itemize}
  \item [(1)] $\Ga$ is the Schl\"{a}fli graph or its complement;
  \item [(2)] $\Ga$ is the Hamming graph $H(s,p^t)$ with $p^t\geq5$ and $st=n$, or $\ov{H(2,p^t)}$ with $p^t\geq5$;
  \item [(3)] $\Ga$ is a normal Cayley graph on $\ZZ_p^n$;
\end{itemize}
\end{proposition}

Schl\"{a}fli graph and Hamming graph are defined in Examples~\ref{PSU} and \ref{Hamming}.

Let $p$ be a prime.
Maru\v{s}i\v{c}~\cite[Theorem 3.4]{Mar} proved that a vertex transitive graph of order $p^2$ or $p^3$ is a Cayley graph, and $2$-arc transitive Cayley graphs of order $p^2$ were classified by Maru\v{s}i\v{c}~\cite[Corollary 2.3]{MP}. This, together with \cite[Corollary 1.3 and 3.5]{Li}, enables us to obtain the following result.

\begin{proposition}\label{2-arc-p^2p^3} For an odd prime $p$, a connected $2$-arc transitive graph of order $p^2$ is $\C_{p^2}$ or $\K_{p^2}$, and a connected $2$-arc transitive graph of order $p^3$ is $\C_{p^3}$, $\K_{p^3}$, or a cover of a complete graph.
\end{proposition}

\section{Examples}\label{Example}

In this section, we give some 2-geodesic transitive graphs, and these contain all $2$-geodesic transitive graphs in question except a new construction in Section~\ref{Newfamily}.

It is well known that $\C_n$, $\K_n$, $\K_{n,n}$, and $\K_{n,n}-n\K_2$ are $2$-geodesic and $2$-arc transitive, and
$\Aut(\C_n)\cong\D_{2n}$, $\Aut(\K_n)\cong\SyG_n$, $\Aut(\K_{n,n})\cong (\SyG_n\times\SyG_n)\rtimes\ZZ_2$ and $\Aut(\K_{n,n}-n\K_2)\cong\SyG_n\times\ZZ_2$. However, the complete multipartite graph $\K_{m[b]}$ with $m\geq 3$ and $b\geq2$, is $2$-geodesic but not $2$-arc transitive, see \cite[Example 2.2(1)]{CJL} and \cite[Lemma 2.4]{Jin17-2}.

\begin{example}\label{multi}
Let $\Ga=\K_{m[b]}$ with $m\geq 3$ and $b\geq2$.
Then $\Aut(\Ga)\cong\SyG_b\wr\SyG_m$ is imprimitive on $V(\Ga)$, and $\Ga$ is $2$-geodesic but not $2$-arc transitive.
\end{example}

Let $G$ be a transitive permutation group on $\Omega$ and let $O$ be an orbit of $G$ on $\Omega\times\Omega$. If $O$ is self-paired, that is, $O=\{(u,v)\ | \ (v,u)\in O\}$, then the {\em orbital graph} of $G$ corresponding to $O$ has vertex set $\Omega$ and edge set $\{\{u,v\}\ |\ (u,v)\in O\}$. The following is the Schl\"{a}fli graph, see \cite[P.26]{Atlas} and \cite[Lemma 3.4]{Jin17}.

\begin{example}\label{PSU}
Let $T\cong\PSU(4,2)$ and $H\cong\ZZ_2^4\rtimes\A_5$. Then $T$ is a primitive permutation group on $[T:H]$, the set of right cosets of $H$ in $T$, and there are two orbital graphs of valency $16$ and $10$ of $T$ on $[T:H]$. Both graphs have order $27$ and are  $2$-geodesic  but not $2$-arc transitive, with automorphism groups isomorphic to $\Aut(T)=\PSU(4,2).\ZZ_2$. The graph with valency $16$ is called Schl\"{a}fli graph and the graph with valency $10$ is the complementary graph of the Schl\"{a}fli graph in $\K_{27}$.
\end{example}

The following is the Hamming graph, see \cite[Section 9.2]{BCN} and \cite[Proposition 2.2]{JDLP}.

\begin{example}\label{Hamming}
Let $d,n$ be positive integers with $d,n\geq2$.
The Hamming graph $H(d,n)$ is defined to be the Cayley graph $\Cay(G,S)$ with $G=\la a_1\ra\times\la a_2\ra\times\cdots\times\la a_d\ra\cong\ZZ_n^d$ and $S=\la a_1\ra^*\cup\la a_2\ra^*\cup\cdots\cup\la a_d\ra^*$. Then $H(d,n)$ is geodesic transitive and  $\Aut(\Ga)\cong\SyG_n\wr\SyG_d$. Furthermore, the complementary graph $\ov{H(2,n)}$ of $H(2,n)$ in the complete graph $\K_{n^2}$, with $n\geq 3$, is $2$-geodesic transitive.
\end{example}

For an odd prime $p$, the extraspecial group $E(p^3)$ of exponent $p$ is defined as follows:
\begin{equation}\label{Ep^3}
E(p^3)=\la a,b,c \mid a^p=b^p=c^p=1,[a,b]=c,[c,a]=[c,b]=1\ra.
\end{equation}

There are two infinite families of $2$-geodesic but not $2$-arc transitive normal Cayley graphs on $E(p^3)$, of which one family is given in the following example which was first constructed in \cite{HFZ}. Another family will be constructed and studied in the next section.

\begin{example}\label{p3-1}
Let $p$ be an odd prime and let $S=\{a^i,b^i\mid i\in\ZZ_p^*\}$. Then $\Cay(E(p^3),S)$ is  normal and $\Aut(\Cay(E(p^3),S))\cong E(p^3)\rtimes ((\ZZ_{p-1}\times \ZZ_{p-1})\rtimes\ZZ_2)$. Furthermore,
$\Cay(E(p^3),S)$ is $2$-geodesic but not $2$-arc transitive.
\end{example}

\proof
Let $\Ga=\Cay(E(p^3),S)$ and $A=\Aut(\Ga)$.
From \cite[Theorem 1.2, Lemmas 3.1 and 3.3]{HFZ}, $\Ga$ is normal and is not $2$-arc transitive with $\Aut(\Ga)\cong E(p^3)\rtimes ((\ZZ_{p-1}\times \ZZ_{p-1})\rtimes\ZZ_2)$. Let $\ZZ_p^*=\langle t\rangle$, and let $\b$ be the automorphism of $E(p^3)$ induced by $a\mapsto a$ and $b\mapsto b^t$. Then $\la\b\ra$ is transitive on $\Ga_2(1)\cap \Ga(a)$, and hence $\Ga$ is $2$-geodesic transitive.\qed

To end this section, we give a proposition which lists all $2$-geodesic transitive graphs of order
$4,8,9,25$ or $27$. This was obtained from the list of all connected symmetric graphs with orders from $2$ to $30$ in \cite{Conder}. For convenience, we write
\[
\begin{array}{l}
\GG_{(27,4)}=\Cay(E(3^3),\{a,b,a^{-1},b^{-1}\}),\\
\GG_{(27,8)}=\Cay(E(3^3), \la b\ra^*\cup_{i\in\ZZ_3}\la b^iab^i\ra^*).
\end{array}
\]
Note that $\GG_{27,4}$ is just the graph for $p=3$ in Example~\ref{p3-1} while $\GG_{(27,8)}$ is just the graph for $p=3$ in Theorem~\ref{NPp3-2}.

\begin{proposition}\label{small-ord}
Let $\Ga$ be a connected $2$-geodesic transitive graph of order $n$.
If $\Ga$ is $2$-arc transitive, then
\begin{itemize}
  \item [(1)] for $n=4$, $\Ga\cong\C_4$ or $\K_4$;
  \item [(2)] for $n=8$, $\Ga\cong\C_8, H(3,2),\K_{4,4}$ or $\K_8$;
  \item [(3)] for $n=9$, $\Ga\cong\C_9$ or $\K_9$;
  \item [(4)] for $n=25$, $\Ga\cong\C_{25}$ or $\K_{25}$;
  \item [(5)] for $n=27$, $\Ga\cong\C_{27}$ or $\K_{27}$.
\end{itemize}
If $\Ga$ is not $2$-arc transitive, then $\Ga$ cannot be of order $4$ and
\begin{itemize}
  \item [(1)] for $n=8$, $\Ga\cong\K_{4[2]}$;
  \item [(2)] for $n=9$, $\Ga\cong H(2,3)$ or $\K_{3[3]}$;
  \item [(3)] for $n=25$, $\Ga\cong H(2,5),\ov{H(2,5)}$, or $\K_{5[5]}$;
  \item [(4)] for $n=27$, $\Ga\cong\GG_{(27,4)},H(3,3),\GG_{(27,8)},\K_{9[3]},\K_{3[9]}$, or the Schl\"{a}fli graph or its complement.
\end{itemize}
\end{proposition}

\section{A family of $2$-geodesic transitive Cayley graphs on $E(p^3)$}\label{Newfamily}

The goal of this section is to prove that the Cayley graph $\Cay(E(p^3), S)$ with $S=\la b\ra^*\cup_{i\in\ZZ_p}\la b^iab^i\ra^*$ is $2$-geodesic transitive but not $2$-arc transitive. To do this, we need the following two lemmas, which will also be used in other Sections.

\begin{lemma}\label{AGL(2,p)}
Let $p\geq5$ be a prime, and let $A=\ZZ_p^2\rtimes \GL(2,p)$. Then  $A\cong\AGL(2,p)$ if and only if every non-identity in the center of $\GL(2,p)$ cannot commute with every element of order $p$ in $\ZZ_p^2$.
\end{lemma}

\proof The necessity is well known. To prove the sufficiency, let $L=\ZZ_p^2$. Then $A=L\GL(2,p)$ and $A/L\cong\GL(2,p)$. To finish the proof, we only need to show that $\soc(A)=L$. Let $N$ be a minimal normal subgroup of $A$. Then $N=T^n$ for a positive integer $n$ and a simple group $T$. Set $Z=Z(\GL(2,p))$. Then $LZ\unlhd A$ and $A/LZ\cong\PGL(2,p)$. Clearly, $\soc(\PGL(2,p))\cong \PSL(2,p)$. By \cite[P. 181, (6.10)]{Hupp}, $\GL(2,p)'=\SL(2,p)\leq \GL(2,p)$.

Suppose that $T$ is a nonabelian simple group. Then $L\cap N=1$ and $LZ\cap N=1$. It follows $N\cong N/(LZ\cap N)\cong NLZ/LZ\unlhd A/LZ\cong\PSL(2,p)$, and then $\PSL(2,p)\cong N\cong N/(N\cap L)\cong NL/L\unlhd A/L\cong \GL(2,p)$, which is impossible.

Thus, $T$ is abelain, and $N=T^n$ with $T\cong \ZZ_q$ for a prime $q$. Suppose that $q\not=p$. Since $L$ is $p$-group, $L\cap \GL(2,p)N=1$, and $\GL(2,p)N\cong \GL(2,p)N/(L\cap \GL(2,p)N)\cong N\GL(2,p)L/L=A/L\cong \GL(2,p)$, forcing $N\GL(2,p)=\GL(2,p)$. Thus, $N\leq \GL(2,p)$ and hence $N\leq Z$. On the other hand, since $N\unlhd A$ and $L\unlhd A$, $N\cap L=1$ implies that $N$ commutes with $L$ pointwise, contradicting the hypothesis. Thus, $q=p$ and $N=\ZZ_p^n$.

If $NL\not=L$, then $NL/L$ is non-trivial $p$-group, and this is impossible because $NL/L\unlhd A/L\cong \GL(2,p)$. Thus, $NL=L$ and $N\leq L$. It follows that every minimal normal subgroup of $A$ is a subgroup of $L$, that is, $\soc(A)\leq L$.

To prove $\soc(A)=L$, the left thing is to show that $L$ is a minimal normal subgroup of $A$. Suppose that $L_1=\la \sigma\ra\leq L$ such that $|L_1|=p$ and $L_1\unlhd A$. Since $\SL(2,p)\leq \GL(2,p)$, we consider the conjugate action of $\SL(2,p)$ on $L_1$.
By \cite[P. 393, (6.3)]{Suzu}, $\SL(2,p)$ has a unique involution, say $\b$, and $\SL(2,p)$ has an element of order $4$, say $\gamma$. Then $\gamma^2=\b$. Clearly, $\b\in Z$. By assumption, $\sigma\not=\sigma^\b\in \la \sigma\ra$, and since $\b$ is an involution, $\sigma^\b=\sigma^{-1}$. Then the orbit set of $\la\b\ra$ on $L_1^*$ is $\Delta=\{\{ \sigma^i,\sigma^{-i} \}\mid i\in\ZZ_{p}^*\}$. Since $\la \b\ra\unlhd \SL(2,p)$, $\SL(2,p)/\la \b\ra\cong\PSL(2,p)$ has an induced natural action on $\Delta$, and then the simplicity of $\PSL(2,p)$ implies the action is trivial because  $\PSL(2,p)$ has no non-trivial action of degree less then $p$ by Proposition~\ref{maxPGL2}. Thus, $\PSL(2,p)$ fixes $\Delta$ pointwise. In particular, $\gamma$ fixes each $\{\sigma^i,\sigma^{-i}\}$ for all $i\in \ZZ_p^*$. It follows that $\gamma^2=\b$ fixes $\la \sigma\ra$ pointwise, contradicting $\sigma^\b=\sigma^{-1}$. This means that $L$ is a minimal normal subgroup of $A$ and hence $\soc(A)=L$, that is,  $A\cong \AGL(2,p)$.\qed


\begin{lemma}\label{autEp^3}
Let $p\geq5$ be a prime. Then we have the following:
\begin{itemize}
  \item[(1)]  Every subgroup $G$ of $\GL(2,p)$ with index less than $p$, contains $\SL(2,p)=\GL(2,p)'$, and for every $m\mid (p-1)$, $\GL(2,p)$ has exactly one subgroup of index $m$.
  \item[(2)] For each subgroup $G$ with $\SL(2,p)\leq G\leq \GL(2,p)$, $\AGL(2,p)$ has exactly
         one conjugate class of subgroups isomorphic to $G$.
     \item[(3)]  $\Aut(E(p^3))\cong \AGL(2,p)$.
\end{itemize}
\end{lemma}

\proof
To prove (1), let $G\leq \GL(2,p)$ such that $|\GL(2,p):G|\bigm|<p$. Then $|\GL(2,p)/Z:GZ/Z|<p$, where  $Z=Z(\GL(2,p))$. Since $\GL(2,p)/Z\cong\PGL(2,p)$, by Proposition~\ref{maxPGL2} we have
$\GL(2,p)/Z=GZ/Z$, implying $\GL(2,p)=GZ$. By \cite[P. 181, (6.10)]{Hupp}, $\GL(2,p)'=\SL(2,p)$, and hence $G\geq G'=(GZ)'=\SL(2,p)$. Now let $|\GL(2,p):G|=m$. Then $m\mid (p-1)$ and $\SL(2,p)\leq G\leq \GL(2,p)$.  Since  $\GL(2,p)/\SL(2,p)\cong \ZZ_{p-1}$, $\GL(2,p)$ has exactly one subgroup of index $m$ for every $m\mid (p-1)$.

\medskip
To prove (2), we may write $\SL(2,p)\leq G\leq \GL(2,p)\leq A=\AGL(2,p)$, and let $B$ be a subgroup of $A$ such that $B\cong G$. We only need to prove that $B$ and $G$ are conjugate in $A$.
Let $L$ be the unique minimal normal subgroup of $A$. Then $L\cong \ZZ_p^2$. Note that $\GL(2,p)'=G'=\SL(2,p)'=\SL(2,p)$. Since $B\cong G$, we have $B'\cong\SL(2,p)$.

Set $M=\SL(2,p)L$. Then $\SL(2,p)\cong M/L\leq A/L\cong \GL(2,p)$. Since $B'\cong \SL(2,p)$, we have $L\cap B'=1$ and $B'L/L\cong B'/(B'\cap L)\cong B'\cong\SL(2,p)$. Then both $M/L$ and $B'L/L$ are subgroups of index $p-1$ in $A/L$. By (1), $M/L=B'L/L$, that is, $M=B'L$.

Denote by $\d$ and $\gamma$ be the center involutions of $\SL(2,p)$ and $B'$, respectively. Then $L\la \d\ra\unlhd M$ and $M/L\la\d\ra\cong \PSL(2,p)$. It follows that $\PSL(2,p)\cong M/L\la\d\ra=B'L/L\la\d\ra\cong B'/(B'\cap L\la\d\ra)$, forcing $B'\cap L\la\d\ra=\la \gamma\ra$. Thus, $\gamma\in L\la \d\ra$. By the Sylow Theorem, $\d$ and $\gamma$ are conjugate in $L\la \d\ra$. Therefore,
the centralizer $C_A(\la \d\ra)$ of $\la\d\ra$ in $A$ is isomorphic to the centralizer $C_A(\la \gamma\ra)$. By Lemma~\ref{AGL(2,p)}, $C_A(\la \d\ra)=\GL(2,p)$, and hence $C_A(\la \gamma\ra)\cong \GL(2,p)$. Thus, $\GL(2,p)$ and $C_A(\la \gamma\ra)$ are conjugate in $A$.

Since $G\cong B$ and $C_A(\la \gamma\ra)\cong \GL(2,p)$, we may write $|\GL(2,p):G|=|C_A(\la \gamma\ra):B|=m$. By (1), $\GL(2,p)$ has a unique subgroup of index $m$, that is, $G$, and $C_A(\la \gamma\ra)$ has a unique subgroup of index $m$, that is, $B$. It follows that $G$ and $B$ are conjugate in $A$.

\medskip
To prove (3), let $A=\Aut(E(p^3))$ and $C=\la c\ra$. By Eq.~(\ref{Ep^3}), $C=Z(E(p^3))$. Then $A$ fixes $C$ and
$E(p^3)\backslash C$ setwise. Let
$x,y\in E(p^3)\backslash C$
with $x\not\in \la y\ra$. Then $\langle x,y\rangle=E(p^3)$. It is easy to see that $x$ and $y$ have the same relations as $a$ and $b$, given in Eq.~(\ref{Ep^3}). By the von Dyck’s Theorem (see \cite[P. 51, 2.2.1]{Rob}), the mapping
\begin{equation}\label{Aut(Ep^3)}
\sigma_{x,y}: a\mapsto x,\ b\mapsto y
\end{equation}
induces an automorphism of $E(p^3)$, denoted by $\sigma_{x,y}$. Clearly, $c^{\sigma_{x,y}}=[a,b]^{\sigma_{x,y}}=[x,y]$. By \cite[Theorem 1.1]{Winter}, $\Aut(E(p^3)\cong \ZZ_p^2\rtimes\GL(2,p)$, and we may assume that $\GL(2,p)\leq A$. \smallskip

\noindent{\bf Claim~1.}\ Let $L=\la \sigma_{a,bc},\sigma_{ac,b}\ra$. Then $A=L\rtimes \GL(2,p)\cong \AGL(2,p)$ and $L$ is the unique minimal normal subgroup of $\Aut(E(p^3))$.\smallskip

It is easy to check that $\sigma_{ac^i,bc^j}\sigma_{ac^k,bc^l}=\sigma_{ac^{i+k},bc^{j+l}}$ for all $i,j,k,l\in\ZZ_p$, and so
$L=\{ \sigma_{ac^i,bc^j}\ |\ i,j\in\ZZ_p\}\cong\ZZ_p^2$. Clearly, $L\leq A$.

Since $A$ fixes $C$ setwise, $A$ has a natural action on the quotient group $E(p^3)/C=\{a^ib^jC\ |\ i,j\in\ZZ_p\}$. Let $K$ be the kernel of this action. Then $K\unlhd A$. Clearly, $L\leq K$. Take $\sigma\in K$. Then $\sigma$ fixes $aC$ and $bC$ setwise, respectively, and hence $\sigma=\sigma_{ac^i,bc^j}$ for some $i,j\in\ZZ_p$, implying $\sigma\in L$. It follows that $K=L$ and $L$ is a normal subgroup of $A$. In particular, $A/L$ has a faithful action on $E(p^3)/C\cong \ZZ_p^2$. Since $L\unlhd A$, we have $L\cap \GL(2,p)=1$ and so $A=L\rtimes \GL(2,p)$.
Then $\GL(2,p)$ induces the automorphism group of the group $E(p^3)/C$, and for convenience, we also write $\GL(2,p)=\Aut(E(p^3)/C)$.

Let $Z$ be the center of $\GL(2,p)$. Take a non-identity $\a\in Z$. Let $q$ be the order of $\a$. Then $q\neq1$, and since $Z\cong\ZZ_{p-1}$, we have $q\mid (p-1)$. Noting that $\a\in Z\leq \GL(2,p)=\Aut(E(p^3)/C)$, $\a$ maps $aC$ to $a^rC$, and $bC$ to $b^rC$, where $r$ is an element of order $q$ in $\ZZ_p^*$. Then $p\nmid r$ and $r\not=1$ as $q\not=1$. It follows that $a^\a=a^rc_1$ and $b^\a=b^rc_2$ for some $c_1,c_2\in C$. In particular, $c^\a=[a,b]^\a=c^{r^2}$. Clearly, $c^r\not=1$ as $p\nmid r$. For every $1\not=\sigma_{ac^s,bc^t}\in L$, we have $c^s\neq1$ or $c^t\not=1$.
If $c^t\not=1$, then $b^{\a\sigma_{ac^s,bc^t}}=b^rc_2c^{tr}\not=b^rc_2c^{tr^2}=b^{\sigma_{ac^s,bc^t}\a}$, and if $c^s\not=1$, then $a^{\a\sigma_{ac^s,bc^t}}=a^rc_1c^{rs}\not=a^rc_1c^{r^2s}=a^{\sigma_{ac^s,bc^t}\a}$. Thus, $\a$ cannot commute with $\sigma_{ac^s,bc^t}$. By Lemma~\ref{AGL(2,p)}, $A\cong \AGL(2,p)$, and hence $L$ is the unique minimal normal subgroup of $\Aut(E(p^3))$. \qed


Now we are ready to prove the main result of this section.

\begin{theorem}\label{NPp3-2}
For an odd prime $p$, let $E(p^3)$ be given as in Eq.~(\ref{Ep^3}), and let $S=\la b\ra^*\cup_{i\in\ZZ_p}\la b^iab^i\ra^*$. Then $\Aut(E(p^3),S)\cong \GL(2,p)$, and $\Cay(E(p^3),S)$ is $2$-geodesic and distance transitive, but not $2$-arc transitive. Furthermore,
\begin{itemize}
  \item[(1)]  For every $m\mid (p-1)$, the unique subgroup of index $m$ in $\Aut(E(p^3),S)$ has $2m+2$ orbits: one orbit has length $1$, one orbit has length $p^2-1$, namely $S$, $m$ orbits have length $(p-1)/m$, and $m$ orbits have length $(p^2-1)(p-1)/m$;
  \item[(2)] $\Aut(E(p^3),S)'\cong \SL(2,p)$ has orbit-set $\{\{c^i\}, c^iS\ |\ i\in\ZZ_p\}$;
     \item[(3)] $|SS|=p^2+(p^2-1)(p-1)$, where $SS=\{s_1s_2\mid s_1,s_2\in S\}$.
\end{itemize}
\end{theorem}

\proof When $p=3$, using {\sc Magma}~\cite{Magma}, one may verify that all results are true. In what follows, we assume that $p\geq 5$.

Let $\Ga=\Cay(E(p^3),S)$ and $C=\la c\ra$. Note that $E(p^3)'=Z(E(p^3))=C$. The following formulae are some basic facts and very useful, which follow from \cite[Theorem 2.1 and Lemma 2.2]{Gore}. For $x,y,z\in E(p^3)$ and $i\in \ZZ_p$, we have $xy=yx[x,y]$, $[x^i,y]=[x,y^i]=[x,y]^i$, $[xy,z]=[x,z][y,z]$, $[x,yz]=[x,y][x,z]$, and for every  integer $k$,  $$(xy)^k=x^ky^k[y,x]^{2^{-1}k(k-1)},$$
where $2^{-1}$ is the inverse of $2$ in $\ZZ_p^*$.

For $i,j\in\ZZ_p$ with $i\not=j$, it is easy to see that $\la b^iab^i\ra\cap \la b^jab^j\ra=1$. Then $|S|=p-1+p(p-1)=p^2-1$. Set
$$T=\la a\ra^*\cup_{i\in\ZZ_p}\la a^iba^i\ra^*.$$ We also have $|T|=p^2-1$.

Let $i\in\ZZ_p^*$. Then $a^iba^i=ba^{2i}c^i$ and $(b^{4^{-1}i^{-1}}ab^{4^{-1}i^{-1}})^{2i}=(b^{2^{-1}i^{-1}}ac^{4^{-1}i^{-1}})^{2i}=
ba^{2i}[a,b^{2^{-1}i^{-1}}]^{2^{-1}2i(2i-1)}c^{2^{-1}}=
ba^{2i}c^{2^{-1}i^{-1}2^{-1}2i(2i-1)}c^{2^{-1}}=ba^{2i}c^i$. It follows that  $a^iba^i=(b^{4^{-1}i^{-1}}ab^{4^{-1}i^{-1}})^{2i}$, implying that $T\subseteq S$. Since $|S|=|T|$,  we have $S=T$.

Let $B=\Aut(E(p^3),S)$. Then $B\leq \Aut(\Ga)_1$ fixes $C$ setwise. Since $S=T$, we have $$S=\la a\ra^*\cup S_1, \mbox{ where } S_1=\cup_{i\in\ZZ_p}\la a^iba^i\ra^*=\la b\ra^*\cup_{i\in\ZZ_p^*}\la b^iab^i\ra^*.$$
Let $x,y\in E(p^3)\backslash C$ with $x\not\in \la y\ra$. Then $\langle x,y\rangle=E(p^3)$. We still use $\sigma_{x,y}$ to denote the automorphism of $E(p^3)$  defined by Eq.~(\ref{Aut(Ep^3)}).

Since $S=\la a\ra^*\cup S_1$, the stabilizer $B_a$ of $a$ in $B$ fixes $S_1$ setwise. Furthermore, for all $i\in\ZZ_p$ we have $\sigma_{a,a^iba^i}\in B_a$ as $S_1=\cup_{i\in\ZZ_p}\la a^iba^i\ra^*$, and for all $j\in\ZZ_p^*$ we have
$\sigma_{a,b^j}\in B_a$ as $S_1= \la b\ra^*\cup_{i\in\ZZ_p^*}\la b^iab^i\ra^*$. It follows that $B_a$ is transitive on $S_1$. Since $E(p^3)$ can be generated by $a$ and $x$, for every $x\in S_1$, $B_a$ is regular on $S_1$. Thus, $|B_a|=p(p-1)$. Similarly, $B_b$ is transitive on  $\cup_{i\in\ZZ_p}\la b^iab^i\ra^*$, and hence $B$ is transitive on $S$. It follow that $|B|=(p^2-1)|B_a|=p(p^2-1)(p-1)=|\GL(2,p)|$.

By Lemma~\ref{autEp^3}, $\Aut(E(p^3))\cong \AGL(2,p)$, and by Claim~1 of the Proof of Lemma~\ref{autEp^3}~(3), $L=\la \sigma_{a,bc}\ra\times \la\sigma_{ac,b}\ra\cong\ZZ_p^2$ is the unique minimal normal subgroup of $\Aut(E(p^3))$. Take $\sigma_{ac^i,bc^j}\in L$ for $i,j\in\ZZ_p$. If $\sigma_{ac^i,bc^j}\in B$, then $\sigma_{ac^i,bc^j}$ fixes $a$ because $S=\la a\ra^*\cup_{i\in\ZZ_p}\la a^iba^i\ra^*$, and $\sigma_{ac^i,bc^j}$ fixes $b$ because $S=\la b\ra^*\cup_{i\in\ZZ_p}\la b^iab^i\ra^*$. It follows that
$\sigma_{ac^i,bc^j}=1$, that is, $L\cap B=1$. This implies that $|LB|=|L||B|=p^3(p^2-1)(p-1)=|\AGL(2,p)|$, and so $\Aut(E(p^3))=L\rtimes B$. In particular, $B\cong \Aut(E(p^3))/L\cong \GL(2,p)$.

Since $S=\la b\ra^*\cup_{i\in\ZZ_p}\la b^iab^i\ra^*$, we have $E(p^3)/C=\{C,sC\ |\ s\in S\}$, and since $|S|=p^2-1$, we have $S\cap sC=\{s\}$ for all $s\in S$. Since $B\cong\GL(2,p)$, we have $\SL(2,p)\cong B'$, and since $B$ fixes $C$, $B$ induces the automorphism group of $E(p^3)/C\cong\ZZ_p^2$. Thus, $B'$ is transitive on
$(E(p^3)/C)^*$, and since $S\cap sC=\{s\}$ for all $s\in S$, $B'$ is transitive on $S$.

Recall that $|B_a|=p(p-1)$ and  $\sigma_{a,b^2}\in B_a$. Then $B_a\cong\ZZ_p\rtimes\ZZ_{p-1}$,  and $B_a$ is transitive on $C\backslash\{1\}$ as $c^{\sigma_{a,b^2}}=[a,b]^{\sigma_{a,b^2}}=c^2$, implying that $B_a$ is transitive on $aC\backslash\{a\}$. It follows that $|B_a:(B_a)_{ac}|=p-1$ and $|(B_a)_{ac}|=p$. On the other hand, $B_{ac}$ fixes $aC\cap S=\{a\}$, and so $|B_{ac}|=|(B_a)_{ac}|=p$. Let $(ac)^{B}$ be the orbit of $B$ containing $ac$. Then $|(ac)^{B}|=|B|/|B_{ac}|=(p^2-1)(p-1)$.

Recall that $B\leq \Aut(\Ga)_1$. Since $1+(p^2-1)+(p^2-1)(p-1)+(p-1)=p^3$, $B$ has exactly four orbits on $V(\Ga)$, that is, $\{1\}$, $S$, $(ac)^{B}$ and $C\backslash\{1\}$. Since $S=\la b\ra^*\cup_{i\in\ZZ_p}\la b^iab^i\ra^*$, we have $aS\cap C\backslash\{1\}=\emptyset$, and since $B$ is transitive on $S$, we have $SS\cap C\backslash\{1\}=\emptyset$. It follows that $\Ga(1)=S$, $\Ga_2(1)=(ac)^B$, and $\Ga_3(1)=C\backslash\{1\}$. In particular, $\Ga$ has diameter $3$ and is distance transitive.
Since $\la a\ra^*\subset S$, $\Ga$ has girth $3$, and hence it is not $2$-arc transitive. Recall that $S=\la a\ra^*\cup S_1$ and $B_a$ is regular on $S_1$. Then $B_a$ is regular on $\{sa\ |\ s\in S_1\}$. Since $a$ is adjacent to all vertices in $\{sa\ |\ s\in S_1\}$, we have $\{as\ |\ s\in S_1\}\not\subseteq S$ (otherwise $\Ga$ is a complete graph). It follows that $\{as\ |\ s\in S_1\}\subseteq \Ga(a)\cap \Ga_2(1)$. Since $\la a\ra^*\subseteq S$, we have $|\Ga(a)\cap \Ga_2(1)|\leq |S_1|$, and hence  $\Ga(a)\cap \Ga_2(1)=\{as\ |\ s\in S_1\}$. Thus, $\Ga$ is $2$-geodesic transitive and $B_a$ is regular on $\Ga(a)\cap \Ga_2(1)$.

To prove part (1), let $m\mid (p-1)$ and let $G$ be the unique subgroup of $B$ of index $m$. Then $G\unlhd B$, and by Lemma~\ref{autEp^3}, $B'\leq G$. Since $B'$ is transitive on $S$, $G$ is transitive on $S$. Since $B_a$ is regular on $\Ga(a)\cap \Ga_2(1)$, $G_a$ is semiregular, and so has $m$ orbits on $\Ga(a)\cap \Ga_2(1)$ as $m=|B:G|=|B_a:G_a|$. Since $B$ is arc-transitive and $G\unlhd B$, $G$ has $m$ orbits on $\Ga_2(1)$ with each length $(p^2-1)(p-1)/m$. Since $C\backslash \{1\}$ is an orbit of $B$ on $V(\Ga)$ and $B\cong \GL(2,p)$, the kernel of $B$ on $C\backslash \{1\}$ is $B'\cong \SL(2,p)$. Then $B/B'$ is regular on $C\backslash \{1\}$, and since $|B/B':G/B'|=|B:G|=m$, $G/B'$ has $m$ orbits on $C\backslash \{1\}$. It follows that $G$ has $m$ orbits on  $C\backslash \{1\}$ with each length $(p-1)/m$.

To prove part (2), take $G=B'=\SL(2,p)$ in the proof of part (1). Then $B'$ has $2p$ orbits: $p$ orbits have length $1$ and $p$ orbits have length $p^2-1$. Furthermore, the $p$
orbits of length $1$ are $\{c^i\}$ for all $i\in\ZZ_p$, and the $p$ orbits of length $p^2-1$ are $c^iS$ for all $i\in \ZZ_p$.

Note that $\la s\ra^*\subset S$ for all $s\in S$. Then $SS=\{1\}\cup \Ga(1)\cup \Ga_2(1)$ and $|SS|=p^2+(p^2-1)(p-1)$, which is part (3). \qed

\section{Proof of Theorem~\ref{prim-p23}}\label{Affine case}


Let $p$ be a prime and let $n\geq 2$ be a positive integer. Write $\Aut(\ZZ_p^n)=\GL(n,p)$,  $Z=Z(\GL(n,p))$, and
$\PGL(n,p)=\GL(n,p)/Z$.
For $z\in\ZZ_p^n$, write $\ov{z}=z^Z$, the orbit of $Z$ containing $z$, and for $T\subseteq \ZZ_p^n$, set $\ov{T}=\{\ov{t}\ |\ t\in T\}$. Recall that for a subgroup $H$ of a group $G$, $H^*$ is the subset of $H$ deleting $1$. Then $\ov{z}=\la z\ra^*$ for $1\not=z\in \ZZ_p^n$, but $\ov{1}=\{1\}$. The kernel of the natural action of $\GL(n,p)$ on $\ov{\ZZ_p^n}$ is $Z$, and $\PGL(n,p)$ is a $2$-transitive permutation group on the point set $\ov{(\ZZ_p^n)^*}=\{\ov{x}\ |\ 1\not=x\in\ZZ_p^n \}$. For $N\leq \GL(n,p)$, write $\overline{N}=NZ/Z$, and then $\overline{N}\leq \PGL(2,p)$. For points $\ov{x_1}, \ov{x_2},\cdots, \ov{x_t}$ in $\ov{(\ZZ_p^n)^*}$, denote by $N_{\ov{x_1}, \ov{x_2},\cdots, \ov{x_t}}$ and $\ov{N}_{\ov{x_1}, \ov{x_2},\cdots, \ov{x_t}}$ the subgroups of $N$ and $\ov{N}$ fixing $\{\ov{x_1}, \ov{x_2},\cdots, \ov{x_t}\}$ pointwise, respectively. First we give some properties on $\Aut(\ZZ_p^3)$, which will be used later.

\begin{lemma}\label{EAp3C}
For a prime $p\geq 5$, let $\ZZ_p^3=\la a\ra\times \la b\ra \times \la c\ra$, $B=\Aut(\ZZ_p^3)=\GL(3,p)$ and $\ov{B}=\PGL(3,p)=\GL(3,p)/Z(\GL(3,p))$. Then we have:
\begin{itemize}
  \item[(1)]  If $\{\ov{a},\ov{b},\ov{c},\ov{u}\}\subseteq \ov{(\ZZ_p^3)^*}$ such that any three points in  $\{\ov{a},\ov{b},\ov{c},\ov{u}\}$ generates $\ZZ_p^3$, then $\ov{B}_{\ov{a},\ov{b},\ov{c}}\cong \ZZ_{p-1}\times \ZZ_{p-1}$ and $\ov{B}_{\ov{a},\ov{b},\ov{c},\ov{u}}=1$;
  \item[(2)] Let $\a\in\ov{B}$ be of order a prime $q$ with $q\mid (p-1)$ and $\b\in\ov{B}$ be of order $p$. Assume that $\a$ commutes with $\b$. Then $\a$ has a fixed point in $\ov{(\ZZ_p^3)^*}$. If further $\a\in\ov{B}_{\ov{a},\ov{b},\ov{c}}$ then $\a$ has exactly
      $3$ or $p+2$ fixed points in $\ov{(\ZZ_p^3)^*}$.
  \item[(3)] $B_{\ov{c}}\cong \AGL(2,p)\times\ZZ_{p-1}$ and $\ov{B}_{\ov{c}}\cong \AGL(2,p)$. Furthermore, $B_{\ov{c}}$ has exactly one conjugate class of subgroups isomorphic to $\SL(2,p)$, of which one has orbit set $\{\{c^i\}, c^i\la a,b\ra^*\ |\ i\in\ZZ_p\}$ in $\ZZ_p^3$.
\end{itemize}
\end{lemma}

\proof For $i,j,k\in\ZZ_p^*$, $\ZZ_p^3$ has an automorphism induced by $a\mapsto a^i$, $b\mapsto b^j$ and $c\mapsto c^k$, denoted by $\sigma_{a^i,b^j,c^k}$. Then $Z=Z(\GL(3,p))=\{\sigma_{a^i,b^i,c^i}\ |\ i\in\ZZ_p^*\}$, $B_{\ov{a},\ov{b},\ov{c}} =\{\sigma_{a^i,b^j,c^k}\ |\ i,j,k\in\ZZ_p^*\}\cong \ZZ_{p-1}\times \ZZ_{p-1}\times \ZZ_{p-1}$, and   $\ov{B}_{\ov{a},\ov{b},\ov{c}}=B_{\ov{a},\ov{b},\ov{c}}/Z\cong
\ZZ_{p-1}\times\ZZ_{p-1}$ as $Z\leq B_{\ov{x},\ov{y},\ov{z}}$. Let $u=a^rb^sc^t$ with $r,s,t\in\ZZ_p^*$. If $\sigma_{a^i,b^j,c^k}\in B_{\ov{a},\ov{b},\ov{c}}$ fixes $\ov{u}$, then
$a^{ir}b^{js}c^{kt}=(a^rb^sc^t)^{\sigma_{a^i,b^j,c^k}}=(a^rb^sc^t)^m$
for some $m\in\ZZ_p^*$, implying $i=j=k=m$ and $\sigma_{a^i,b^j,c^k}\in Z$. It follows that $\ov{B}_{\ov{a},\ov{b},\ov{c},\ov{u}}=1$. This completes the proof of part~(1).

To prove part (2), we suppose to the contrary that $\a$ has no fixed point in $\ov{(\ZZ_p^3)^*}$. Since $|\ov{(\ZZ_p^3)^*}|=p^2+p+1$, we have $q \mid (p^2+p+1)$, and hence $q$ is a divisor of $(p^2+p+1)-(p-1)=p^2+2$.
Since $(p^2+2)-p(p-1)=p+2$ and $p+2-(p-1)=3$, we have $q=3$. Since $p\nmid (p^2+p+1)$, $\b$ has a fixed point in $\ov{(\ZZ_p^3)^*}$, say $\ov{u}$. Let $\ov{u}^{\la \a\ra}=\{\ov{u},\ov{v},\ov{w}\}$.
Since $\a$ has no fixed point, $\la \ov{u},\ov{v},\ov{w} \ra=\ZZ_p^2$ or $\ZZ_p^3$. If $\la \ov{u},\ov{v},\ov{w} \ra=\ZZ_p^2$, then $3\mid (p+1)$, forcing $3$ is a divisor of $(p+1)-(p-1)=2$, a contradiction. Thus, $\la \ov{u},\ov{v},\ov{w} \ra=\ZZ_p^3$.
Since $\a$ commutes with $\b$, $\b$ fixes $\ov{u}$, $\ov{v}$ and  $\ov{w}$, that is, $\b\in\ov{B}_{\ov{u},\ov{v},\ov{w}}$, contradicting part~(1).

Let $\a\in\ov{B}_{\ov{a},\ov{b},\ov{c}}$.  Then $\a$ fixes $\ov{a},\ov{b}$ and $\ov{c}$, and hence $\a=\sigma_{a^i,b^j,c^k}Z$ for some $i,j,k\in\ZZ_p^*$. Since $\sigma_{a^{i^{-1}},b^{i^{-1}},c^{i^{-1}}}\in Z$, we have $\a=\sigma_{a^i,b^j,c^k}Z=\sigma_{a,b^r,c^s}Z$, where $r=i^{-1}j\in\ZZ_p^*$ and $s=i^{-1}k\in\ZZ_p^*$. Note that $o(\a)=q$. If $r=1$ then $s\not=1$ and $\a$ has the fixed-point set $\ov{\la a,b\ra^*}\cup \{\ov{c}\}$, and if $s=1$ then $r\not=1$ and $\a$ has the fixed-point set $\ov{\la a,c\ra^*}\cup \{\ov{b}\}$. For $r\not=1$ and $s\not=1$, if $r=s$ then $\a$ has the fixed-point set $\{\ov{a}\}\cup \ov{\la b,c\ra^*}$, and if $r\not=s$ then
$\a$ has the fixed-point set $\{\ov{a},\ov{b},\ov{c}\}$. Thus, $\a$ has exactly $3$ or $p+2$ fixed points in $\ov{(\ZZ_p^3)^*}$.

To prove part~(3), let $\ZZ_p^3=\la x,y,z\ra$ and denote by $\delta_{x,y,z}$ the automorphism of $\ZZ_p^3$ induced by $a\mapsto x$, $b\mapsto y$ and $c\mapsto z$. Let $L$ be the kernel of $B_{\ov{c}}$ on ${\ov{c}}$. Then $|B_{\ov{c}}/L|\leq p-1$, and since $Z$ is transitive on ${\ov{c}}$, we have $B_{\ov{c}}=L\times Z$ and $\ov{B}_{\ov{c}}\cong L$.

Write $\Aut(\langle a,b\ra)=\GL(2,p)$ and $\GL(2,p)'=\SL(2,p)$. Set $G=\{\delta_{a^\a,b^\a,c}\ |\ \a\in \GL(2,p)\}$ and $H=\{\delta_{a^\a,b^\a,c}\ |\ \a\in \SL(2,p)\}$. Then $H\leq G\leq L$ with $H\cong\SL(2,p)$ and $G\cong\GL(2,p)$. Clearly, $G$ and $H$ have the same orbit set on $\ZZ_p^3$, that is, $\{\{c^i\},c^i\la a,b\ra^*\ |\ i\in\ZZ_p\}$.

Let $M\cong\SL(2,p)$ be a subgroup of $B_{\ov{c}}$. Then $M=M'\leq B_{\ov{c}}'=L'\leq L$. Note that Lemma~\ref{autEp^3}~(2) means that  $\AGL(2,p)$ has one conjugate class of subgroups isomorphic to $\SL(2,p)$. Thus, to finish the proof of
part~(3), it suffices to show that $L\cong \AGL(2,p)$.

Let $C=\la c\ra$. Since $L$ fixes $C$ pointwise, $L$ has a natural action on $\ZZ_p^3/C$, and let $K$ be the kernel of the action of $L$ on the quotient group $\ZZ_p^3/C$. Then $K$ fixes $aC$ and $bC$, and hence $K=\{ \delta_{ac_1,bc_2,c}\ |\ c_1,c_2\in C\}\cong\ZZ_p^2$.
Since $K\unlhd L$, we have $G\cap K=1$, and since $|L/K|\leq |\Aut(\ZZ_p^3/C)|=|G|$, we have $L=K\rtimes G$.

Note that $Z(G)=\{\delta_{a^r,b^r,c}\ |\ r\in\ZZ_p^*\}$. Take $1\not=\delta_{ac_1,bc_2,c}\in K$ and $1\not=\delta_{a^r,b^r,c}\in Z(G)$. Then $r\not=1$, and either $c_1\not=1$ or $c_2\not=1$. If  $c_1\not=1$ then  $a^{\delta_{ac_1,bc_2,c}\delta_{a^r,b^r,c}}=a^rc_1\not=a^r(c_1)^r=a^{\delta_{a^r,b^r,c}\delta_{ac_1,bc_2,c}}$, and if $c_2\not=1$ then $b^{\delta_{ac_1,bc_2,c}\delta_{a^r,b^r,c}}=b^rc_2\not=b^r(c_2)^r=b^{\delta_{a^r,b^r,c}\delta_{ac_1,bc_2,c}}$. Thus, $\delta_{ac_1,bc_2,c}\delta_{a^r,b^r,c}\neq\delta_{a^r,b^r,c}\delta_{ac_1,bc_2,c}$.
By Lemma~\ref{AGL(2,p)}, $L\cong \AGL(2,p)$, as required.  \qed

We also need some properties on cliques of certain normal Cayley graphs on $\ZZ_p^n$.

\begin{lemma}\label{clique}
Let $\Ga=\Cay(\ZZ_p^n,S)$ be a connected $2$-geodesic but not $2$-arc transitive normal Cayley graph on $\ZZ_p^n$. Then
\begin{itemize}
  \item[(1)] $\Aut(G,S)$ contains the center of $\Aut(\ZZ_p^n)$, and for every positive integer $m$, if $x\in \Ga_m(1)$ then $\la x\ra^*\subseteq \Ga_m(1)$;
  \item[(2)] Let $K_1$ and $K_2$ be subgroups of $\ZZ_p^n$. Assume that the induced subgroups $[K_1]$ and $[K_2]$ of $\Ga$ are complete graphs and every vertex in $K_1$ is adjacent to every vertex in $K_2$. Then the induced subgraph $[K_1K_2]$ of $\Ga$ is a complete graph;

    \item[(3)] $G$ has a subgroup $H$ such that $H^*\subseteq S$ and $[H]$ is a maximal clique of $\Ga$.
\end{itemize}
\end{lemma}

\proof For $r\in\ZZ_p^*$, let $\sigma_r$ be the automorphism of $\ZZ_p^n$ which maps every element $x\in\ZZ_p^n$ to $x^r$. Then the center $Z(\Aut(\ZZ_p^n))=\{\sigma_r\ |\ r\in\ZZ_p^*\}$. By Proposition~\ref{nor-2-geo}, $\la s\ra^*\subset S=\Ga_1(1)$ for each $s\in S$, and hence $\sigma_r$ fixes $S$ setwise. Thus, $Z(\Aut(\ZZ_p^n))\leq \Aut(\ZZ_p^n,S)\leq \Aut(\Ga)$, which implies that for every positive integer $m$, if $x\in \Ga_m(1)$ then $\la x\ra^*\subseteq \Ga_m(1)$.

To prove part (2), assume that $[K_1]$ and $[K_2]$ are complete graphs and every vertex in $K_1$ is adjacent to every vertex in $K_2$. Clearly, $K_1K_2\leq \ZZ_p^n$.  Take two distinct vertices $x_1y_1,x_2y_2\in K_1K_2$ with $x_1,x_2\in K_1$ and $y_1,y_2\in K_2$.
Then $x_1x_2^{-1}\not=y_1^{-1}y_2$. Since $x_1x_2^{-1}\in K_1$ is adjacent to $y_1^{-1}y_2\in K_2$, we have $x_1x_2^{-1}(y_1^{-1}y_2)^{-1}\in S$, that is, $x_1y_1(x_2y_2)^{-1}\in S$. Thus, $x_1y_1$ is adjacent to $x_2y_2$, and so $[K_1K_2]$ is a complete graph.

To prove part~(3), assume that $H$ is a maximal subgroup of $G$ such that $[H]$ is a complete graph. Then $H^*\subseteq S$. If $[H]$ is not a maximal clique of $\Ga$, then there is $x\not\in H$ such that $x$ is adjacent to all elements of $H$. In particular, $x\in S$ as $x$ is adjacent to $1$. By (1), $\la x\ra^*\subseteq S=\Ga_1(1)$ and hence $\la x\ra$ is a complete graph. Since $Z(\Aut(\ZZ_p^n))\leq \Aut(\Ga)$, every vertex in $\la x\ra$ is adjacent to every vertex in $H$, and by (2), $[\la x\ra H]$ is a complete graph, contradicting the maximality of $H$. It follows that $[H]$ is a maximal clique of $\Ga$. \qed

In the next lemma, we consider $2$-geodesic but not $2$-arc transitive normal Cayley graphs on $\ZZ_p^2$.

\begin{lemma}\label{HA-p^2}
Let $p\geq 5$ be a prime and let $\Ga=\Cay(\ZZ_p^2,S)$ be a connected $2$-geodesic but not $2$-arc transitive normal Cayley graph. Then $\Aut(\Ga)$ cannot be primitive on $V(\Ga)$.
\end{lemma}

\proof  Let $B=\GL(2,p)=\Aut(\ZZ_p^2)$ and $A=\Aut(\Ga)$ . Since $\Ga$ is normal, Propositions~\ref{Aut-Cay} implies $A=R(\ZZ_p^2)\rtimes\Aut(\ZZ_p^2,S)$ and $A_1=\Aut(\ZZ_p^2,S)\leq B$, where $A_1$ is the stabilizer of $1$ in $A$.
In particular, $\Ga\not=\K_{p^2}$.
By Lemma~\ref{clique}~(1), $Z=Z(\GL(2,p))\leq A_1$. Denote by $\overline{\Ga}$ the quotient graph of $\Ga$ corresponding to $Z$, that is, the graph with the orbits of $Z$ as vertices and with two orbits adjacent if and only if there is an edge of $\Ga$ between the two orbits. Then $V(\overline{\Ga})=\ov{\ZZ_p^2}$. Since $Z\leq A_1\leq B$, the kernel of the action of $A_1$ on $\ov{\Ga}$ is $Z$. Then $\overline{A_1}\leq \Aut(\overline{\Ga})_{\ov{1}}$, and $\overline{A_1}$ is a permutation group on $\ov{(\ZZ_p^2)^*}$. By Lemma~\ref{clique}~(1), $\overline{\Ga}_m(\ov{1})=\ov{\Ga_m(1)}$ for every $m\geq 1$.

Suppose to the contrary that $A$ is primitive on $V(\Ga)$. Since $\Ga$ is connected, $\la S\ra=\ZZ_p^2$, and so we may let $a,b\in S$ with $\ZZ_p^2=\la a,b\ra$. Then $\ov{a}\subseteq S$ and $\ov{b}\subseteq S$. Since $\ZZ_p^2=\langle a\rangle\langle b\rangle\subseteq SS$, $\Ga$ has diameter $2$.
Note that $|\ov{(\ZZ_p^2)^*}|=p+1$. Suppose $|\ov{S}|=p+1$. By  Lemma~\ref{clique}~(1), $S=(\ZZ_p^2)^*$ and $\Ga$ is a complete graph, a contradiction. Suppose $|\ov{S}|=p$. Then $|\ov{\Ga}_2(\ov{1})|=1$ and hence $\langle \Ga_2(1)\rangle\cong\ZZ_p$. Therefore, $A_1$ fixes a subgroup of order $p$ in $\ZZ_p^2$, contrary to the primitivity of $A$ on $V(\Ga)$. Suppose $|\ov{S}|=2$. Then $S=\ov{a}\cup\ov{b}$ and $\Ga=H(2,p)$. By Example~\ref{Hamming}, $\Aut(\Ga)\cong (\SyG_p\times\SyG_p)\rtimes \ZZ_2$, which is impossible because $\Aut(\Ga)$ has no normal subgroup isomorphic to $\ZZ_p^2$. Suppose $|\ov{S}|=p-1$. Then $|\ov{\Ga}_2(\ov{1})|=(p+1)-(p-1)=2$ and $\Ga=\ov{H(2,p)}$, the complement of $H(2,p)$ in the complete graph $\K_{p^2}$, which is also impossible because $\Aut(\ov{H(2,p)})\cong \Aut(H(2,p))$. It follows that $3\leq |\ov{S}|\leq p-2$.

Clearly,  $\ov{(\ZZ_p^2)^*}=\{\ov{a},\ov{ba^i}\ |\ i\in\ZZ_p\}$. Since $\Ga$ has diameter $2$, $\ov{\Ga}_2(\ov{1})=\{\ov{ba^i}\ |\ \ov{ba^i}\not\in \ov{S},i\in\ZZ_p^*\}$, and for $\ov{ba^i}\in \ov{\Ga}_2(\ov{1})$, we have $\ov{ba^i}\in \ov{\Ga}_1(\ov{b})$ because $b$ is adjacent $ba^i$. Thus, $$\ov{\Ga}_2(\ov{1})\subseteq \ov{\Ga}_1(\ov{b}).$$
Let $\ov{u},\ov{v}\in \ov{\Ga}_2(\ov{1})$. Then $\ov{u},\ov{v}\in \ov{\Ga}_1(\ov{b})$, and since $Z\leq \Aut(\Ga)$, $b$ is adjacent to some vertex $u^i\in \ov{u}$ and $v^j\in \ov{v}$. This implies that  $(1,b,u^i)$ and $(1,b,v^j)$ are two 2-geodesics of $\Ga$, and since $\Ga$ is $2$-geodesic transitive, the stabilizer $A_{1b}$ of $b$ in $A_1$
has an element mapping $(1,b,u^i)$ to $(1,b,v^j)$. It follows that $\ov{A_{1b}}$ is transitive on $\ov{\Ga}_2(\ov{1})$, and hence  $(\ov{A_1})_{\ov{b}}$ is transitive on $\ov{\Ga}_2(\ov{1})$. Since $\Ga$ is arc-transitive, $A_1$ is transitive on $\Ga_1(1)=S$, and so $\ov{A_1}$ is transitive on $\ov{S}$.
Thus, $|\ov{A_1}|=|\ov{S}||(\ov{A_1})_{\ov{b}}|$, and the transitivity of $(\ov{A_1})_{\ov{b}}$ on $\ov{\Ga}_2(\ov{1})$ implies that $|\ov{S}||\ov{\Ga}_2(\ov{1})|\mid |\ov{A_1}|$. Since $\Ga$ has diameter $2$,  $|\ov{\Ga}_2(\ov{1})|=(p+1)-|\ov{S}|$, and therefore, $\ell(p+1-\ell) \mid |\ov{A_1}|$, where $\ell=|\ov{S}|$. Since $3\leq \ell\leq p-2$, we have $3\leq p+1-\ell\leq p-2$. Thus, $(\ell(p+1-\ell),p)=1$ and hence $\ell(p+1-\ell)\mid |\ov{A_1}|_{p'}$. Note that the function $f(t)=t(p+1-t)$ has the minimal value $3(p-2)$ when $3\leq t\leq p-2$.
It follows that $$3\leq \ell\leq p-2,\ 3(p-2)\leq \ell(p+1-\ell) \mbox{ and } \ \ell(p+1-\ell) \mid |\ov{A_1}|_{p'}.$$

Suppose $\PSL(2,p)\leq \ov{A_1}$. By Proposition~\ref{maxPGL2}, every transitive action of $\ov{A_1}$ has degree at least $p$, and so $|\ov{S}|=\ell\geq p$, a contradiction. Thus, $\ov{A_1}\leq M$, where $\PSL(2,p)\nleq M$ and $M$ is a maximal subgroup of $\PSL(2,p)$ or $\PGL(2,p)$.
By Proposition~\ref{maxPGL2}, $M\cong \ZZ_p\rtimes\ZZ_{p-1},\ZZ_p\rtimes\ZZ_{(p-1)/2}, \D_{2(p-1)}, \D_{p-1}, \D_{2(p+1)}, \D_{p+1},\A_4, \SyG_4$ or $\A_5$,
and since $\ell(p+1-\ell) \mid |\ov{A_1}|_{p'}$, $\ell(p+1-\ell)$ is divisor of $2(p+1)$, $2(p-1)$, $24$ or $60$.

Suppose $\ell(p+1-\ell)\mid 2(p-1)$. Since $\ell(p+1-\ell)\geq 3(p-2)$ and $3(p-2)=2(p-1)+p-4$, we have $p\leq 4$, contradicting the hypothesis $p\geq 5$.

Suppose $\ell(p+1-\ell)\mid 24$. Since $\ell(p+1-\ell)\geq 3(p-2)$,
we have either $p=5$ and $\ell=3$, or $p=7$ and $3\leq\ell\leq5$, of which each is impossible.

Suppose $\ell(p+1-\ell)\mid 60$. Then $\ov{A_1}\leq \A_5$. If $\ov{A_1}<\A_5$, then $\ov{A_1}$ is a subgroup of a maximal subgroup of $\A_5$, and by Atlas~\cite{Atlas},
$\ell(p+1-\ell)$ is a divisor of $12$ or $20$,
which are impossible by a similar argument to the above case $\ell(p+1-\ell)\mid 24$. Thus, $\ov{A_1}=\A_5$, and since $\ov{S}$ and $\ov{\Ga}_2(\ov{1})$ are orbits of $\ov{A_1}$, we have $\ell=|\ov{S}|\geq 5$ and $p+1-\ell=|\ov{\Ga}_2(\ov{1})|\geq 5$. It follows that $p\geq 9$, and since $3(p-2)\leq \ell(p+1-\ell)\leq 60 $, we have $p\leq 22$. Thus, $p=11,13,17$ or $19$, and for each case, $\ell(p+1-\ell)\nmid 60$, a contradiction.

Suppose $\ell(p+1-\ell)\mid 2(p+1)$. Then $\ell(p+1-\ell)\geq 3(p+1-3)=2(p+1)+p-8$. If $p\geq 9$ then $\ell(p+1-\ell)$ is not a divisor of $2(p+1)$. Thus, $p=5$ or $7$, and $\ell(p+1-\ell)\mid 2(p+1)$ implies that $p=7$ and $\ell=4$. It follows $\ov{A_1}= M\cong \D_{16}$, and $|\ov{\Ga}_1(\ov{1})|=|\ov{\Ga}_2(\ov{1})|=4$. Since $\ov{\Ga}_1(\ov{1})$ and $\ov{\Ga}_2(\ov{1})$ are orbits of $\ov{A_1}$,
the unique involution of $\ov{A_1}$ fixes every point in  $\ov{(\ZZ_7^2)^*}$,
contradicting that $\ov{A_1}$ is  a permutation group on $\ov{(\ZZ_7^2)^*}$. \qed

Now we consider $2$-geodesic but not $2$-arc transitive normal Cayley graphs on $\ZZ_p^3$.

\begin{lemma}\label{HA-p3}
Let $p\geq 5$ be a prime and let $\Ga=\Cay(\ZZ_p^3,S)$ be a connected $2$-geodesic but not $2$-arc transitive normal Cayley graph. Then $\Aut(\Ga)$ cannot be primitive on $V(\Ga)$.
\end{lemma}

\proof Let $B=\Aut(\ZZ_p^3)=\GL(3,p)$, $Z=Z(\GL(3,p))$ and $\ov{B}=\PGL(2,p)$. Let $A=\Aut(\Ga)$. By Proposition~\ref{Aut-Cay}, $A=R(\ZZ_p^3)\rtimes\Aut(\ZZ_p^3,S)$, and $A_1=\Aut(\ZZ_p^3,S)\leq B$, implying $\Ga\not=\K_{p^3}$. By Lemma~\ref{clique}~(1), $Z\leq A_1$. Let $\overline{\Ga}$ be the quotient graph of $\Ga$ corresponding to $Z$. Then $V(\overline{\Ga})=\ov{\ZZ_p^3}$, and since $Z\leq A_1\leq B$, we have $\overline{A_1}\leq \Aut(\overline{\Ga})_{\ov{1}}$. Furthermore, $\overline{A_1}\leq \ov{B}=\PGL(3,p)$ and $\overline{\Ga}_m(\ov{1})=\ov{\Ga_m(1)}$ for every $m\geq 1$.

Suppose, by way of contradiction, that $A$ is primitive on $V(\Ga)$. Let $\O=\{\O_1,\cdots,\O_t\}$ be the set of orbits of $\ov{A_1}$ on $\ov{(\ZZ_p^3)^*}$. The primitivity of $A$ implies that $\la \O_i\ra=\ZZ_p^3$ for each $1\leq i\leq t$, forcing $|\ov{\O_i}|\geq 3$ for each $1\leq i\leq t$. In particular, $|\ov{\Ga}_m(\ov{1})|\geq 3$ and  $\la \ov{\Ga}_m(\ov{1})\ra=\ZZ_p^3$, for every $1\leq m\leq \diam(\Ga)$.

\medskip
Denote by $O_p(\ov{A_1})$  the maximal normal $p$-subgroup of $\ov{A_1}$.

\medskip
\noindent{\bf Claim:} $O_p(\ov{A_1})=1$.

Suppose $O_p(\ov{A_1})\not=1$.
If $O_p(\ov{A_1})$ has no fixed point in $\ov{(\ZZ_p^3)^*}$, then $p$ is a divisor of $|\ov{(\ZZ_p^3)^*}|=p^2+p+1$, which is impossible. Now assume that $O_p(\ov{A_1})$ has a fixed point in $\O_j$ for some $1\leq j\leq m$. Since $O_p(\ov{A_1})\unlhd \ov{A_1}$, $O_p(\ov{A_1})$ fixes $\O_j$ pointwise. Since $\la \O_j\ra =\ZZ_p^3$, we may assume that $\ov{x},\ov{y},\ov{z}\in\O_j$ and $\ZZ_p^3=\la x,y,z\ra$. Then  $O_p(\ov{A_1})$ fixes $\ov{x},\ov{y}$ and $\ov{z}$, and hence  $O_p(\ov{A_1})\leq \ov{B}_{\ov{x},\ov{y},\ov{z}}$, contradicting Lemma~\ref{EAp3C}~(1). Thus, $O_p(\ov{A_1})=1$, as claimed.

\medskip

Since $\Ga$ is connected, $\la S\ra=\ZZ_p^3$, and so we may assume $\ZZ_p^3=\la a,b,c\ra$ with $\{\ov{a}, \ov{b}, \ov{c}\}\subseteq \ov{S}$.
If $S=\ov{a}\cup\ov{b}\cup\ov{c}$,
then $\Ga=H(3,p)$, and by Example~\ref{Hamming}, $\Aut(\Ga)\cong (\SyG_p\times\SyG_p\times\SyG_p)\rtimes \SyG_3$, which is impossible because $\Aut(\Ga)$ has no normal subgroup isomorphic to $\ZZ_p^3$. Thus, we may assume that $\ell:=|\ov{S}|\geq 4$. By Lemma~\ref{clique}~(3), $\ZZ_p^3$ has a subgroup $H$ such that $H^*\subseteq S$ and $[H]$ is a maximal clique of $\Ga$. Since $\Ga$ is not a complete graph, $H\cong \ZZ_p^2$ or $\ZZ_p$.

\medskip
\noindent{\bf Case 1}: $H\cong\ZZ_p^2$.

Without loss of any generality, we may let $H=\la a,b\ra$. Then $\ov{H^*}\subseteq \ov{S}$ and $\ov{c}\in \ov{S}$, implying  $|\ov{S}|\geq p+2$, and since $\ZZ_p^3=H\langle c\rangle\subseteq SS$, $\Ga$ has diameter $2$.
Note that $|\ov{(\ZZ_p^3)^*}|=p^2+p+1$
and $\ov{(\ZZ_p^3)^*}=\ov{H^*}\cup \{\ov{ch}\ |\ h\in H\}$, where $|\ov{H^*}|=p+1$ and $|\{\ov{ch}\ |\ h\in H\}|=p^2$. Then $\ov{\Ga}_2(\ov{1})=\{\ov{ch}\ |\ 1\not=h\in H, \ov{ch}\not\in \ov{S}\}$. Since $c$ is adjacent $ch$, we have $$\ov{\Ga}_2(\ov{1})\subseteq \ov{\Ga}_1(\ov{c}).$$

Let $\ov{u},\ov{v}\in \ov{\Ga}_2(\ov{1})$. Then $\ov{u},\ov{v}\in \ov{\Ga}_1(\ov{c})$, and since $Z\leq \Aut(\Ga)$, $c$ is adjacent to some vertex $u^i\in \ov{u}$ and $v^j\in \ov{v}$. Since $\Ga$ is $2$-geodesic transitive, the stabilizer $A_{1c}$ of $c$ in $A_1$
has an element mapping the $2$-geodesic $(1,c,u^i)$ to $(1,c,v^j)$. Thus, $\ov{A_{1c}}$ is transitive on $\ov{\Ga}_2(\ov{1})$, and hence  $(\ov{A_1})_{\ov{c}}$ is transitive on $\ov{\Ga}_2(\ov{1})$. The arc-transitivity of $\Ga$ implies that $A_1$ is transitive on $\Ga_1(1)=S$, and so $\ov{A_1}$ is transitive on $\ov{S}$. It follows that
$|\ov{S}||\ov{\Ga}_2(\ov{1})|\mid |\ov{A_1}|$. Since $\Ga$ has diameter $2$, we have $|\ov{S}|+|\ov{\Ga}_2(\ov{1})|=p^2+p+1$.
Then $\ell(p^2+p+1-\ell) \mid |\ov{A_1}|$, where $p+2\leq \ell=|\ov{S}|\leq p^2+p$. If $\ell\geq p^2+p+1-2$, then $|\ov{\Ga}_2(\ov{1})|=p^2+p+1-\ell\leq 2$, contradicting $|\ov{\Ga}_2(\ov{1})|\geq 3$. If $\ell=p^2+p+1-3$, then $|\ov{\Ga}_2(\ov{1})|=3$ and $\Ga=\ov{H(3,p)}$, which is impossible because $\Aut(\Ga)\cong (\SyG_p\times\SyG_p\times\SyG_p)\rtimes \SyG_3$ has no normal subgroup isomorphic to $\ZZ_p^3$. Therefore,
$$7\leq p+2\leq \ell\leq p^2+p-3,\ \mbox{ where } \ell=|\ov{S}|.$$
For every $m\geq 1$, it is easy to see that $$f(\ell,p):=\ell (p^2+p+1-\ell)\geq m(p^2+p+1-m), \ \mbox { when } m\leq \ell \leq p^2+p+1-m.$$ In particular, $f(\ell,p)\geq 4(p^2+p-3)$. Note that $f(\ell,p)\mid |\ov{A_1}|$.

By Proposition~\ref{maxPGL3}~(2), $\PSL(3,p)$ has no non-trivial transitive action with degree less than least $p^2+p+1$. This implies that if $\PSL(3,p)\leq \ov{A_1}$, then
$|\ov{S}|=\ell=p^2+p+1$, a contradiction.
Thus, $\ov{A_1}\leq M$, where $\PSL(3,p)\nleq M$ and $M$ is a maximal subgroup of $\PGL(3,p)$ or $\PSL(3,p)$. By Claim, $O_p(\ov{A_1})=1$, and by Proposition~\ref{maxPGL3}, $|\ov{A_1}|$ is a divisor of  $p(p+1)(p-1)^2$, $6(p-1)^2$, $3(p^2+p+1)$, $3^3\cdot 8$, $168$ or $360$, where the last three cases correspond to $M=(\ZZ_3^2:\Q_8).\ZZ_3$ with $p\equiv 1(\mod 3)$, $M=\PSL(2,7)$ with $p\equiv 1,2,4(\mod 7)$, or $M=\A_6$ with $p\equiv 1,4(\mod 15)$, respectively.

Suppose $|\ov{A_1}|\mid 6(p-1)^2$. Then $f(\ell,p)\mid 6(p-1)^2$. If  $7\leq \ell \leq p^2+p-6$ then
$f(\ell,p)\geq 7(p^2+p-6)>6(p-1)^2$, which is impossible. Since $7\leq \ell\leq p^2+p-3$, we have $\ell=p^2+p-5,p^2+p-4$ or $p^2+p-3$, implying $f(\ell,p)\nmid 6(p-1)^2$, a contradiction.

Suppose $|\ov{A_1}|\mid 3(p^2+p+1)$. Then $f(\ell,p)\mid 3(p^2+p+1)$, which is impossible because
$f(\ell,p)\geq 4(p^2+p-3)>3(p^2+p+1)$.

Suppose $|\ov{A_1}|$ is a divisor of $3^3\cdot 8$, $168$ or $360$.
Then $f(\ell,p)\mid 3^3\cdot 8$, $168$ or $360$, of which each is impossible for $p\geq 11$, because otherwise $f(\ell,p)\geq 4(p^2+p-3)\geq 516$. Thus, $p=7$ or $5$. For $3^3\cdot 8$,
$M=(\ZZ_3^2:\Q_8).\ZZ_3$
with $p\equiv 1(\mod 3)$, and for $168$ or $360$, $M=\PSL(2,7)$ with $p\equiv 1,2,4(\mod 7)$ or
$M=\A_6$ with $p\equiv 1,4(\mod 15)$,
of which the latter two are impossible because $p\geq 11$. Thus,
$M=(\ZZ_3^2:\Q_8).\ZZ_3$ and $p=7$, which is also impossible because $f(\ell,7)=\ell(57-\ell)$ is not a divisor of $3^3\cdot 8$ for every $7\leq \ell\leq 53$.

Suppose $|\ov{A_1}|\mid p(p+1)(p-1)^2$. Then $f(\ell,p)\mid p(p+1)(p-1)^2$. By Proposition~\ref{maxPGL3},
$M=\ZZ_p^2:\GL(2,p)$, $\ZZ_p^2:\frac{1}{d}\GL(2,p)$ or $\PGL(2,p)$.
Since $O_p(\ov{A_1})=1$, we have $\ov{A_1}=\ov{A_1}/(\ov{A_1}\cap \ZZ_p^2)\cong \ov{A_1}\ZZ_p^2/\ZZ_p^2\leq M/\ZZ_p^2\cong \GL(2,p)$, $\frac{1}{d}\GL(2,p)$, or $\PGL(2,p)$. Thus, we may assume  $\ov{A_1}\leq N$, where $N=\GL(2,p)$, $\frac{1}{d}\GL(2,p)$, or $\PGL(2,p)$. Then $|Z(N)|\mid (p-1)$ and $\ov{A_1}/(\ov{A_1}\cap Z(N))\cong \ov{A_1}Z(N)/Z(N)\leq \PGL(2,p)$.

Assume $p\nmid |\ov{A_1}|$. Then $\ov{A_1}/(\ov{A_1}\cap Z(N))\not= \PGL(2,p)$. By Proposition~\ref{maxPGL2}, $|\ov{A_1}|$ and so $f(\ell,p)$ is a divisor of $2(p-1)^2$, $2(p^2-1)$, $24(p-1)$ or $60(p-1)$, of which the last one corresponds to
$\ZZ_{p-1}.\A_5$ with $p\equiv \pm1(\mod 10)$. The former two are impossible because $f(\ell,p)\geq 4(p^2+p-3)> 2(p^2-1)$, and the latter two are impossible for $p\geq 17$ because $f(\ell,p)\geq 4(p^2+p-3)>4(p+1)(p-1)>72(p-1)$.
Furthermore, $24(p-1)$ cannot happen because $f(\ell,p)>4(p+1)(p-1)\geq 24(p-1)$, and for $60(p-1)$, we have $p=11$ or $13$ as
$p\equiv \pm1(\mod 10)$. If $7\leq \ell\leq (p^2+p-6)$ then $f(\ell,p)\geq 7(11^2+11-6)=882>60(13-1)>60(11-1)$, which is impossible. It follows that $p=11$ or $13$ with $\ell=p^2+p-5,p^2+p-4$ or $p^2+p-3$, which are also impossible because otherwise $f(\ell,p)\nmid 60(p-1)$.

Assume $p\mid |\ov{A_1}|$. If $Z(N)\cap \ov{A_1}$ has an element $\a$ of order a prime $q$, then $q\mid (p-1)$, and by Lemma~\ref{EAp3C}~(2), $\a$ has a fixed point in $\ov{(\ZZ_p^3)^*}$, which is contained in $\ov{S}$ or $\ov{\Ga}_2(\ov{1})$.
Since $\la\a\ra\unlhd \ov{A_1}$, $\a$ fixes $\ov{S}$ or $\ov{\Ga}_2(\ov{1})$ pointwise, and so $\a$ has fixed-point set $\ov{S}$ or $\ov{\Ga}_2(\ov{1})$. Since $\la \ov{S}\ra=\la\ov{\Ga}_2(\ov{1})\ra=\ZZ_p^3$, $\a$ fixe $\ov{x},\ov{y},\ov{z}$ such that $\ZZ_p^3=\la x,y,x\ra$. Since $7\leq \ell\leq p^2+p-3$, again by Lemma~\ref{EAp3C}~(2) we have either $|\ov{S}|=p+2$ or $|\ov{\Ga}_2(\ov{1})|=p+2$. If follows that $f(\ell,p)=(p+2)(p^2-1)$ is a divisor of $ p(p+1)(p-1)^2$, which is clearly impossible.
Thus, $\ov{A_1}\cap Z(N)=1$, and hence $\ov{A_1}\cong
\ov{A_1}Z(N)/Z(N)\leq \PGL(2,p)$.

Thus, $|\ov{A_1}|\mid p(p+1)(p-1)$ and $f(\ell,p)\mid p(p+1)(p-1)$. If $p+1\leq \ell\leq p^2+ p+1-(p+1)=p^2$ then $f(\ell,p)\geq p^2(p+1)>p(p+1)(p-1)$, which is impossible. It follows that $7\leq \ell\leq p$ or $p^2+1\leq \ell\leq p^2+p-2$, which implies that either $\ell=p,p^2+1$, or $(f(\ell,p),p)=1$.
If $\ell=p$ or $p^2+1$ then $f(\ell,p)=p(p^2+1)\mid p(p+1)(p-1)$, which is impossible. If $(f(\ell,p),p)=1$ then $\ell(p^2+p+1-\ell)\mid (p+1)(p-1)$, which is also impossible because $\ell(p^2+p+1-\ell)\geq 4(p^2+p-3)>(p+1)(p-1)$.

\medskip
\noindent{\bf Case 2}: $H\cong\ZZ_p$.

Recall that $\ov{a}\cup \ov{b}\cup \ov{c}\subseteq S$.
Since $\ZZ_p^3=\la a\ra\la b\ra\la c\ra\leq SSS$, $\Ga$ has diameter at most $3$. Take $\ov{x},\ov{y}\in \ov{S}$ with $\ov{x}\not=\ov{y}$. If $\ov{x}$ is adjacent to every vertex in $\ov{y}$, then Lemma~\ref{clique} implies that $\la x,y\ra$ is a maximal clique of order $p^2$, a contradiction.
Thus, if $\ov{c}\not=\ov{x}$ and $\ov{c}\not=\ov{y}$, then $c$ is not adjacent to some $x^i\in \ov{x}$ and $y^j\in\ov{y}$, and since $\Ga$ is $2$-geodesic transitive, $A_{1c}$ has an element mapping the $2$-geodesic $(c,1,x^i)$ to $(c,1,y^j)$. Thus, $(\ov{A_1})_{\ov{c}}$ is transitive on $\ov{S}\backslash\{\ov{c}\}$, and the transitivity of $\ov{A_1}$ on $\ov{S}$ means that $\ov{A_1}$ is $2$-transitive on $\ov{S}$.

Similarly, if $x$ is not adjacent to $y^i$ and $y^j$ in $\ov{y}$, then
$(x,1,y^i)$ and $(x,1,y^j)$ are two $2$-geodesics, and hence
$A_{1x}$ is transitive on $\ov{y}_{x-}$, where
$$\ov{y}_{x-}=\{y^k\ |\ y^k\in \ov{y},\ x \mbox{ is not adjacent to } y^k\}.$$
Since $A_1\leq B=\Aut(\ZZ_p^3)$, the block stabilizer $(A_{1x})_{\ov{y}}$ of $\ov{y}$ in $A_{1x}$ is transitive on $\ov{y}_{x-}$. Then the restriction $(A_{1x})_{\ov{y}}^{\ov{y}}$ of $(A_{1x})_{\ov{y}}$ on $\ov{y}$ is a subgroup of $\ZZ_{p-1}$, and hence regular on $\ov{y}_{x-}$, implying $|(A_{1x})_{\ov{y}}^{\ov{y}}|=|\ov{y}_{x-}|$. By the $2$-transitivity of $\ov{A_1}$ on $\ov{S}$, we may assume that
$$|\ov{y}_{x-}|=r\ \mbox{ for all } x,y\in S \mbox{ with } \ov{x}\not=\ov{y}, \ \mbox{ where } r\mid (p-1).$$
Denote by $\ov{y}_{x+}$ the set of vertices in $\ov{y}$ that are adjacent to $x$ in $\Ga$. Then $$|\ov{y}_{x+}|=(p-1)-r=(p-1)(1-\frac{r}{p-1}).$$

Recall that $\ov{A_1}$ is $2$-transitive on $\ov{S}$ and $\ov{A_1}\leq \PGL(3,p)$. Let $K$ be the kernel of $\ov{A_1}$ on $\ov{S}$. Since $\la S\ra\cong\ZZ_p^3$, by Lemma~\ref{EAp3C}~(1) we have $K\leq \ZZ_{p-1}\times\ZZ_{p-1}$.

Suppose $O_p(\ov{A_1}/K)\not=1$. Let $O_p(\ov{A_1}/K)=O/K$. Then $O\unlhd \ov{A_1}$. Since $K\leq \ZZ_{p-1}\times\ZZ_{p-1}$, we have $O=KO_p$, where $O_p\not=1$ is a Sylow $p$-subgroup of $O$. We prove that $O_p$ has trivial conjugate action on $K$ by induction on $|K|$. Let $q\mid |K|$ be a prime. Then $q\mid (p-1)$, and the subgroup $L$ generated by all elements of order $q$ in $K$ is characteristic in $K$. Since $K\unlhd O$, we have $L\unlhd O$, and $L\cong\ZZ_q$ or $\ZZ_q^2$. Since $p\geq 5$, we have $p\nmid |\Aut(L)|$, and hence $O_p$ has a trivial conjugate action on $L$. Then $O_p$ has a natural action on $K/L$, and by inductive hypothesis, we may assume that the action is trivial. Since $(|O_p|,|L|)=1$, $O_p$ has a trivial conjugate action on $K$
(see \cite[P. 131, 18.6]{Hupp}),
that is, $O=K\times O_p$. Then $O_p$ is characteristic in $O$, and since $O\unlhd \ov{A_1}$, we have $O_p\unlhd \ov{A_1}$, contradicting Claim. Thus, $O_p(\ov{A_1}/K)=1$.

If $\PSL(3,p)\leq \ov{A_1}/K$, then $|\ov{S}|=p^2+p+1$ because $\PSL(3,p)$ has the minimum transitive degree $p^2+p+1$ by Proposition~\ref{maxPGL3}, contradicting that
$\Ga\not=\K_{p^3}$. Since $\ell=|\ov{S}|\geq 4$ and $O_p(\ov{A_1}/K)=1$,  by Proposition~\ref{2-transPGL3} we have the following lists:
\begin{itemize}
  \item[List A:] $(\ell, \ov{A_1}/K)=(4,\A_4)$, $(4,\SyG_4)$ or $(9,L)$ with $L\leq (\ZZ_3^2:\Q_8).\ZZ_3$ and $p \equiv 1(\mod 3)$;
  \item[List B:] $(\ell, \ov{A_1}/K)=(5,\A_5)$,
    $(6,\A_6)$, $(10,\A_6)$, $(8,\PSL(2,7))$, $(11,\PSL(2,11))$,
      $(p+1,L)$ with $L=\PSL(2,p)$, or $\PGL(2,p)$, and if $X/K=\A_6$, $\PSL(2,7)$ or $\PSL(2,11)$ then $p\equiv 1,4(\mod 15)$, $p\equiv 1,2,4 (\mod 7)$ or $p=11$, respectively.
\end{itemize}
The $2$-transitivity of $\ov{A_1}$ on $\ov{S}$ implies that the induced subgraph $[\ov{S}]$ of $\ov{S}$ in $\ov{\Ga}$ is a complete graph or an empty graph.

Assume that $[\ov{S}]$ is a complete graph. Then $r\not=p-1$, and for $\ov{x},\ov{y}\in \ov{S}$ with $\ov{x}\not=\ov{y}$ we have   $|\ov{y}_{x+}|=(p-1)(1-r/(p-1))\geq (p-1)/2$ as $r\mid (p-1)$. Thus, $|\ov{\la x,y\ra}\cap \ov{S}|\geq (p-1)/2+2$. Since $\{\ov{a},\ov{b},\ov{c}\}\subseteq \ov{S}$, we have $|\ov{\la a,b\ra}\cap \ov{S}|\geq (p-1)/2+2$, $|\ov{\la a,c\ra}\cap \ov{S}|\geq (p-1)/2+2$ and $|\ov{\la b,c\ra}\cap \ov{S}|\geq (p-1)/2+2$. Since $\la a,b\ra\cap \la a,c\ra=\la a\ra$,
$\la a,b\ra\cap \la b,c\ra=\la b\ra$ and $\la a,c\ra\cap \la b,c\ra=\la c\ra$, we have $\ell=|\ov{S}|\geq 3(p-1)/2+6-3=3(p+1)/2$.
In particular, $\ell\geq 9$ as $p\geq 5$, and $\ell>p+1$. By Lists A and B, either $(\ell, \ov{A_1}/K)=(9,L)$ with $L\leq (\ZZ_3^2:\Q_8).\ZZ_3$ and  $p \equiv 1(\mod 3)$, or $(\ell, \ov{A_1}/K)=(11,\PSL(2,11))$ with $p=11$, of  which both are impossible because otherwise $p\geq 7$ and $\ell\geq 3(p+1)/2\geq 12$.

Assume that $[\ov{S}]$ is an empty graph. Then $r=p-1$, and for $\ov{x},\ov{y}\in \ov{S}$ with $\ov{x}\not=\ov{y}$, we have $|\ov{y}_{x-}|=p-1$ and so $\ov{\la x,y\ra}\cap \ov{\Ga}_2(\ov{1})=p-1$.
Let $\ov{c}\not=\ov{x}$ and $\ov{c}\not=\ov{y}$. Then $|\ov{\la c,x\ra}\cap\ov{\Ga}_2(\ov{1})|=|\ov{\la c,y\ra}\cap \ov{\Ga}_2(\ov{1})|=p-1$, and
since $\la c,x\ra\cap\la c,y\ra=\la c\ra$, we have $|(\ov{\la c,x\ra}\cup\ov{\la c,y\ra})\cap\ov{\Ga}_2(\ov{1})|=2(p-1)$. Note that $(\ov{\la c,x\ra}\cup\ov{\la c,y\ra})\cap\ov{\Ga}_2(\ov{1})\subseteq \ov{\Ga}_1(\ov{c})$.  Since $\ell=|\ov{S}|$, the arbitrariness of $\ov{x}$ and $\ov{y}$ implies that $|\ov{\Ga}_1(\ov{c})\backslash\{\ov{1}\}|\geq (\ell-1)(p-1)$. Conversely, for $\ov{z}\in \ov{\Ga}_1(\ov{c})\backslash\{\ov{1}\}$, $c$ is adjacent to $z^i$ with $i\in\ZZ_p^*$, implying $\ov{cz^{-i}}\in \ov{S}$ and $\ov{z}\in \ov{\la c,cz^{-i}\ra}$. Noting that $|\ov{\la a,b\ra}\cap \ov{\Ga}_2(\ov{1})|=p-1$ and $\ov{\la a,b\ra}\cap \ov{\Ga}_1(\ov{c})=\{\ov{1}\}$, we have
$$|\ov{\Ga}_1(\ov{c})\backslash\{\ov{1}\}|=(\ell-1)(p-1) \mbox{ and } |\ov{\Ga}_2(\ov{1})|\geq
(\ell-1)(p-1)+(p-1)
=\ell(p-1).$$
Note that $\ov{\Ga}_1(\ov{c})\backslash\{\ov{1}\}\subseteq \ov{\Ga}_2(\ov{1})$.  Then $(\ov{A_1})_{\ov{c}}$ is transitive  $\ov{\Ga}_1(\ov{c})\backslash\{\ov{1}\}$. It follows that  $$(\ell-1)(p-1)\mid |(\ov{A_1})_{\ov{c}}|\ \mbox{ and }\ell(\ell-1)(p-1)\mid |\ov{A_1}|.$$

Since $\ell=|\ov{S}|\geq 4$, we may let $\ov{u}\in \ov{S}$ with $\ov{u}\not\in\{\ov{a},\ov{b},\ov{c}\}$. If $u\in \la a,b\ra$ then $\ov{u}\in \ov{\la a,b\ra}$, which is impossible because $|\ov{\la a,b\ra}\cap \ov{\Ga}_2(\ov{1})|=p-1$. Then $u\not\in \la a,b\ra$, and similarly, $u\not\in \la a,c\ra$ and $u\not\in \la b,c\ra$. By Lemma~\ref{EAp3C}~(1), $\ov{B}_{\ov{u},\ov{a},\ov{b},\ov{c}}=1$ and so $K=1$. Thus, Lists A and B hold for $K=1$.
Note that $p\geq5$ and $|\ov{A_1}|\geq \ell(\ell-1)(p-1)$.
If $\ell=4$ then $|\ov{A_1}|\geq 3\cdot 4^2>|\SyG_4|$; if $\ell=5$ then $|\ov{A_1}|\geq 80>|\A_5|$; if $\ell=6$ or $10$ then $p\equiv 1,4(\mod 15)$ and $|\ov{A_1}|\geq 30\cdot 18>|\A_6|$; if $\ell=8$ then $p\equiv 1,2,4 (\mod 7)$ and $|\ov{A_1}|\geq 560>|\PSL(2,7)$; if $\ell=9$ then $p\equiv 1(\mod 3)$ and $|\ov{A_1}|\geq 9\cdot 8\cdot 6>|(\ZZ_3^2:\Q_8).\ZZ_3|$; and if $\ell=11$ then $p=11$ and  $|\ov{A_1}|\geq 11\cdot 10^2>|\PSL(2,11)|$.
By Lists A and B, we have $(\ell,\ov{A_1})=(p+1,\PSL(2,p)$) or $(p+1,\PGL(2,p))$. Then $p^2\nmid |\ov{A_1}|$ and  $|\ov{\Ga}_2(\ov{1})|\geq \ell(p-1)=p^2-1$. Since $\ov{S}=\ell=p+1$ and $|\ov{(\ZZ_p^3)^*}|=p^2+p+1$, we have $|\ov{\Ga}_2(\ov{1})|=p^2$ or $p^2-1$: for the former, $p^2\mid |\ov{A_1}|$, a contradiction, and for the latter, $\ov{\Ga}_3(\ov{1})=1$, contradicting $\ov{\Ga}_3(\ov{1})\geq 3$. \qed

\noindent{\bf Proof of Theorem~\ref{prim-p23}.}\ Let $\Ga$ be a 2-geodesic transitive but not $2$-arc transitive Cayley graph on $\ZZ_p^n$ with $p$ a prime and $n\leq 3$. By \cite{JLW}, every $2$-geodesic transitive graph of order $p$ is isomorphic to $\C_p$ or $\K_p$, which are also $2$-arc transitive. So there do not exist 2-geodesic transitive but not 2-arc transitive graphs of order $p$.

In what follows, assume that $n\geq 2$. If $p\geq 5$, then by Lemmas~\ref{HA-p^2} and ~\ref{HA-p3}, $\Aut(\Ga)$ is imprimitive on $V(\Ga)$. If $p\leq 3$, then $p^n=4,8,9$ or $27$ and then by Proposition~\ref{small-ord}, $\Ga$ is isomorphic to $H(2,3)$, $\K_{3[3]}$, $\K_{4[2]}$, $\GG_{(27,4)}$, $H(3,3)$, $\GG_{(27,8)}$, $\K_{9[3]}$, $\K_{3[9]}$, the Schl\"{a}fli graph or its complement. Among these ten graphs, by using {\sc Magma}~\cite{Magma}, it is easy to verify that $\Ga$ is a Cayley graph on $\ZZ_p^2$ or $\ZZ_p^3$ and $\Aut(\Ga)$ is primitive on $V(\Ga)$ if and only if $\Ga\cong H(2,3)$ or $H(3,3)$. This completes the proof.\qed

\section{Proofs of Theorems~\ref{Thm-p^2} and \ref{Thm-p^3}}\label{MainProof}

In this section, we prove Theorems~\ref{Thm-p^2} and \ref{Thm-p^3}. First we prove the following lemma.

\begin{lemma}\label{cover-Kp} For an odd prime $p$, let
$\Ga$ be a connected $2$-geodesic but not $2$-arc transitive graph of order $p^2$ or $p^3$, and let $A=\Aut(\Ga)$. Let $1\neq N\lhd A$ be intransitive on $V(\Ga)$ and let $\Ga$ be a cover of $\Ga_N$. Then $\Ga_N{~\not \cong~}\K_p$.
\end{lemma}

\proof Suppose, by way of contradiction, that $\Ga_N\cong\K_p$. Then  $\Ga$ has valency $p-1$, and $|N|=p$ or $p^2$ as $|V(\Ga)|=p^2$ or $p^3$. Note that the complete multipartite graphs $\K_{p[p]}$, $\K_{p^2[p]}$ and $\K_{p[p^2]}$ cannot be a cover of $\K_p$, because non of the valencies of them is $p-1$. Since $p$ is odd, only (2) or (4) of Proposition~\ref{redu} can occur. For (2), $\Ga$ is bipartite, and since $\Ga$ is $2$-geodesic transitive, it is $2$-arc transitive, a contradiction. For (4), $N$ is semiregular, and $\Ga_N$ is $A/N$-arc transitive but not $(A/N,2)$-arc transitive.
Since $\Ga_N\cong\K_p$, $A/N$ is 2-transitive but not 3-transitive on $V(\Ga_N)$, and by \cite[P. 99]{Dixon}, $A/N$ is one of the following groups:
\begin{itemize}
  \item[(a)] $A/N=\AGL(1,p)\cong\ZZ_p\rtimes\ZZ_{p-1}$;
  \item[(b)] $A/N= \A_p$ or $\Sy_p$;
  \item[(c)] $A/N= \M_{11}$ and $p=11$, or $\M_{23}$ and $p=23$;
  \item[(d)] $A/N= \PSL(2,11)$ and $p=11$;
  \item[(e)] $\PSL(d,q)\leq A/N\leq \P\Gamma L(d,q)$, where $d\geq2$ and $\frac{q^d-1}{q-1}=p$.
\end{itemize}

Assume that case (a) occurs. Then $A=N.\AGL(1,p)$.
Let $G$ be a Sylow $p$-subgroup of $A$. Then $N\leq G$ and $G/N$ is the unique Sylow $p$-subgroup of $\AGL(1,p)$. In particular, $|G|=p^2$ or $p^3$,  $G\lhd A$, and $G$ is regular on $V(\Ga)$. Thus, $\Ga$ is a normal Cayley graph on $G$, say $\Cay(G,S)$.
By Proposition~\ref{nor-2-geo}, $\langle s\ra^*\subseteq S$ for all $s\in S$ and $s$ has order $p$. Since $\Ga$ has valency $p-1$, we have $S=\la s\ra^*$ for some $s\in S$, and by the connectivity of $\Ga$, $G=\langle S\ra=\langle s\rangle\cong\ZZ_p$, contradicting that $|G|=p^2$ or $p^3$.

Cases (b) and (c) cannot occur, because otherwise $A/N$ is  $3$-transitive on $V(\Ga_N)$.

Assume that case (d) or (e) occurs. Then $T:=\soc(A/N)=\PSL(d,q)$.
Clearly, $p^2\nmid |T|$.
Furthermore, $T=X/N$ for some $X\unlhd A$, and hence $X=N.T\unlhd A$. Since $|N|=p$ or $p^2$, by Proposition~\ref{mult} we have that $X$ is a central extension, $X=NX'$ and $X'\cong M.T$ with $M\leq\Mult(T)\cap N$.

Now we claim $M=1$. Since $M\leq N$, $M$ is a $p$-group.
For case (d), $A/N=\PSL(2,11)$, $p=11$, and $T\cong\PSL(2,11)$. By \cite[P. 302, Table 4.1]{Gore-82}, $\Mult(\PSL(2,11))\cong\ZZ_2$, and since $M$ is a $p$-group, $M=1$, as claimed. For case (e), $T\cong\PSL(d,q)$.
If $d=2$, then $\frac{q^d-1}{q-1}=p$ implies that $q$ is a power of $2$ and $T$ is $3$-transitive on $V(\Ga_N)$, forcing that $A/N$ is $3$-transitive on $V(\Ga_N)$, a contradiction. Thus, $d\geq3$, and $\frac{q^d-1}{q-1}=p$ implies that $d$ is an odd prime and $d<p$. If $(d,q)\neq (3,2)$ and $(3,4)$, by \cite[P. 302, Table 4.1]{Gore-82} we have  $\Mult(\PSL(d,q))\cong\ZZ_{(d,q-1)}$, and since $M$ is a $p$-group, $M=1$, as claimed.
If $(d,q)=(3,4)$, then $\frac{q^d-1}{q-1}$ is not a prime, a contradiction.
If $(d,q)=(3,2)$, then $p=7$, and again by \cite[P. 302, Table 4.1]{Gore-82}, $\Mult(\PSL(3,2))\cong\ZZ_2$, forcing $M=1$, as claimed.

It follows that $M=1$ and $X'\cong T$. Since $X'$ is characteristic in $X$, we have $X'\unlhd A$. Since $p^2\nmid |T|$, $X'$ is intransitive on $V(\Ga)$, and since $p$ is odd and $\Ga$ has valency $p-1$, Proposition~\ref{redu} implies that $X'\cong T$ is semiregular on $V(\Ga)$, contradicting that $|V(\Ga)|=p^2$ or $p^3$. \qed

\vskip0.05in
{\noindent\bf Proof of Theorem~\ref{Thm-p^2}.}
By Examples in Section~\ref{Example}, all graphs in Theorem~\ref{Thm-p^2} are 2-geodesic transitive of order $p^2$.

Assume that $\Ga$ is a connected $2$-geodesic transitive graph of order $p^2$ and let $A=\Aut(\Ga)$.
If $\Ga$ is 2-arc transitive, Proposition~\ref{2-arc-p^2p^3} implies that $\Ga\cong\C_{p^2}$ or $\K_{p^2}$, which is  Theorem~\ref{Thm-p^2}~(1). Suppose that $\Ga$ is not 2-arc transitive. If $p\leq5$, then Proposition~\ref{small-ord} implies that $\Ga\cong H(2,3),H(2,5),\ov{H(2,5)},\K_{3[3]}$ or $\K_{5[5]}$. Suppose that $p\geq7$.

Assume that $A$ is quasiprimitive on $V(\Ga)$. By Proposition~\ref{quas-2-geo}, either $\Ga\cong H(2,p)$ or $\ov{H(2,p)}$ that is listed in Theorem~\ref{Thm-p^2}~(2); or $A$ is primitive on $V(\Ga)$ and $\Ga=\Cay(\ZZ_p^2,S)$ is a normal Cayley graph, which is impossible by Lemma~\ref{HA-p^2}.

Assume that $\Ga$ is not quasiprimitive on $V(\Ga)$.
Then $A$ has a nontrivial intransitive normal subgroup $N$ on $V(\Ga)$.
By Proposition~\ref{redu}, either $\Ga\cong\K_{p[p]}$ with $p\geq3$, or $N$ is semiregular on $V(\Ga)$ and $\Ga$ is a cover of $\Ga_N$. For the former, $\Ga$ is listed in Theorem~\ref{Thm-p^2}~(2). For the latter, $|V(\Ga_N)|=p$ and $\Ga_N$ is $(A/N,2)$-geodesic transitive. Then $\Ga_N$ is circulant, and by Proposition~\ref{2-geo-cir}, $\Ga_N\cong\K_p$, which is impossible by  Lemma~\ref{cover-Kp}. \qed

To prove Theorem~\ref{Thm-p^3}, we first classify $2$-arc transitive graphs of order prime-cube.

\begin{lemma}\label{2arc}
Let $p$ be a prime and let $\Ga$ be a connected $2$-arc transitive graph of order $p^3$. Then $\Ga$ is isomorphic to $H(3,2)$, $\K_{4,4}$, $\C_{p^3}$
 or $\K_{p^3}$.
\end{lemma}

\proof If $p=2$ or 3, then Proposition~\ref{small-ord} implies that  $\Ga\cong\C_8,H(3,2),\K_{4,4},\K_{8},\C_{27}$ or $\K_{27}$. Let $p\geq5$.
By Proposition~\ref{2-arc-p^2p^3}, either $\Ga\cong\C_{p^3}$ or $\K_{p^3}$, or
$\Ga$ is a cover of a complete graph. Let $A=\Aut(\Ga)$. If $\Ga$ is a cover of a complete graph, then $A$ has a nontrivial normal subgroup $N$ such that $N$ is semiregular on $V(\Ga)$ and $\Ga$ is a $2$-arc transitive cover of $\Ga_N$ with $\Ga_N =\K_p$ or $\K_{p^2}$. Clearly, $N\cong\ZZ_p$, $\ZZ_{p^2}$ or $\ZZ_p\times\ZZ_p$. By \cite[Theorem 1.1]{DMO}, $p=2$, contradicting that $p\geq5$. It follows that $\Ga\cong\C_{p^3}$ or $\K_{p^3}$.
\qed

\vskip0.05in
{\noindent\bf Proof of Theorem~\ref{Thm-p^3}.}
By Examples in Section~\ref{Example} and Theorem~\ref{NPp3-2}, all graphs in Theorem~\ref{Thm-p^3} are $2$-geodesic transitive of order $p^3$.

Let $\Ga$ be a $2$-geodesic transitive graph of order $p^3$ and let $A=\Aut(\Ga)$. If $\Ga$ is $2$-arc transitive, by Lemma~\ref{2arc}, Theorem~\ref{Thm-p^3}~(1) holds. In what follows we assume that $\Ga$ is not $2$-arc transitive. If $p=2$ or $3$, Proposition~\ref{small-ord} implies that $\Ga\cong \K_{4[2]}$, $\GG_{(27,4)}$, $H(3,3)$, $\GG_{(27,8)}$, $\K_{3[9]}$, $\K_{9[3]}$, or the Schl\"{a}fli graph or its complement. Thus, we let $p\geq5$.

Assume that $A$ is quasiprimitive on $V(\Ga)$. By Proposition~\ref{quas-2-geo}, either $\Ga$ is the Hamming graph $H(3,p)$ that is Theorem~\ref{Thm-p^3}~(2)-(iii), or $A$ is primitive on $V(\Ga)$ and $\Ga$ is a normal Cayley graph on $\ZZ_p^3$, which is impossible by Lemma~\ref{HA-p3}.

Assume that $A$ is not quasiprimitive on $V(\Ga)$.
Then $A$ has a nontrivial maximal intransitive normal subgroup $N$ and so $A/N$ is quasiprimitive on $V(\Ga_N)$. By Proposition~\ref{redu}, either $\Ga\cong\K_{p^2[p]}$ or $\K_{p[p^2]}$ that is Theorem~\ref{Thm-p^3}~(2)-(ii), or $N$ is semiregular on $V(\Ga)$ and $\Ga$ is a cover of $\Ga_N$. To finish the proof, assume that $\Ga$ is a cover of $\Ga_N$, and we only need to show that Theorem~\ref{Thm-p^3}~(2)~(iv) occurs. Thus, $N$ is semiregular on $V(\Ga)$, and $\Ga$ and $\Ga_N$ have the same valency.

Note that $|V(\Ga_N)|=p$ or $p^2$, and $\Ga_N$ is $(A/N,2)$-geodesic transitive. If $|V(\Ga_N)|=p$, then $\Ga_N$ is circulant, and by Proposition~\ref{2-geo-cir}, $\Ga_N\cong \C_{p}$ or $\K_p$, of which the latter is impossible by  Lemma~\ref{cover-Kp}. For the former, $\Ga\cong \C_{p^3}$, contradicting that $\Ga$ is not $2$-arc-transitive.
Thus, we may let $|V(\Ga_N)|=p^2$ and $N\cong\ZZ_p$.

Let $M/N$ be a minimal normal subgroup of $A/N$. Then $M\unlhd A$ and $M/N=T_1\times T_2\times\cdots\times T_t$, where $T_i\cong T$ is a simple group. Since $A/N$ is quasiprimitive on $V(\Ga_N)$, $M/N$ is transitive on $V(\Ga_N)$, and since $\Ga_N$ is $2$-geodesic transitive and $|V(\Ga_N)|=p^2$, $\Ga_N$ is listed in Theorem~\ref{Thm-p^2}.
By Proposition~\ref{cover}, $\Ga_N{~\not\cong~}\K_{p[p]}$.
Since $\Ga$ is 2-geodesic but not 2-arc transitive, $\Ga$ has girth 3, and hence  $\Ga_N$ has girth $3$, which implies $\Ga_N{~\not\cong~}\C_{p^2}$. It follows from Theorem~\ref{Thm-p^2} that
\begin{equation}\label{Eq3}
 \Ga_N\cong\K_{p^2}, H(2,p) \text{~or~} \ov{H(2,p)}.
\end{equation}

\medskip
\noindent{\bf Claim 1:} $M/N$ is abelian.

Suppose that $M/N$ is non-abelian and we will obtain a contradiction. In this case, $T$ is a non-abelian simple group and $M/N=T_1\times T_2\times\cdots\times T_t$ with $T_i\cong T$ for all $i$. We first prove $t=1$, that is, $M/N$ is non-abelian simple.

Suppose $t\geq 2$. Recall that $M/N$ is transitive on $V(\Ga_N)$ and $|V(\Ga_N)|=p^2$. If some $T_i$ is transitive on $V(\Ga_N)$, say $T_1$, then a stabilizer $(T_2)_\a$ with $\a\in V(\Ga_N)$ fixes every vertex in $V(\Ga_N)$, because $T_1$ commutes with $T_2$ pointwise. Thus, $T_2$ is semiregular on $V(\Ga_N)$ and hence $|T_2|$ is a divisor of $p^2=|V(\Ga_N)|$, contradicting that $T_1$ is non-abelian simple. It follows that each orbit of $T_i$ with $1\leq i\leq p$, has length $p$. Let $$\Delta=\{\Delta_1,\cdots,\Delta_p\}$$ be the orbit set of $T_1$ on $V(\Ga_N)$. Then $|\Delta_j|=p$ for all $1\leq j\leq p$, and $\Delta$ is a complete imprimitive block system of $M/N$.

Since $M/N$ is transitive on $V(\Ga_N)$, some $T_i$ is transitive on $\Delta$, say $T_2$. Let $\a_1\in \Delta_1$. Since $T_1$ commutes with $T_2$ pointwise, the transitivity of $T_2$ on $\Delta$ implies that $(T_1)_{\a_1}$ fixes a vertex $\a_i$ in $\Delta_i$ for all $1\leq i\leq p$. Thus, $(T_1)_{\a_1}=(T_1)_{\a_2}=\cdots=(T_1)_{\a_p}$.

Consider the natural action of $M/N$ on $\Delta$ and let $K$ be the kernel of the action. Clearly, $T_1\leq K$ and $T_2\nleq K$. If $T_1\not=K$, then some $T_i$ with $i\geq 3$ is a subgroup of $K$ because $K\unlhd M/N$, says $T_3\leq K$. Then $\Delta$ is also the orbit set of $T_3$. Since $T_1$ and $T_3$ commutes pointwise, $(T_1)_{\a_1}$ fixes every vertex in $\Delta_1$, and similarly, $(T_1)_{\a_i}$ fixes every vertex in $\Delta_i$ for each $1\leq i\leq p$. Since $(T_1)_{\a_1}=(T_1)_{\a_2}=\cdots=(T_1)_{\a_p}$, $T_1$ is semiregular on $V(\Ga_N)$, a contradiction. Thus, $K=T_1$.

It follows that $(M/N)/T_1$ has a faithful action on $\Delta$, and since $|\Delta|=p$, $(M/N)/T_1$ is  isomorphic to a subgroup of the symmetric group $S_p$. Thus, $p^2 \nmid (M/N)/T_1$. It follows that $t=2$ and so $M/N=T_1\times T_2$ with  $p^2 \nmid |T_2|$. On the other hand, the simplicity of $T_1$ implies that $T_1$ is isomorphic to a subgroup of the symmetric group on $\Delta_1$, and similarly, $p^2 \nmid |T_1|$. Thus, $p^3\nmid |T_1\times T_2|$.

Since $M/N=T_1\times T_2$, we may let $X/N=T_2$ with $N\leq X\unlhd M$.
Since $T_2$ is transitive on $\Delta$, it has a subgroup of index $p$. By Proposition~\ref{simple-ind}, we have
\begin{equation*}
T_2\cong\M_{11},\M_{23},\A_p,\PSL(2,11)\text{~with~} p=11,\text{~or~}\PSL(d,q)\text{~with~}(q^d-1)/(q-1)=p.
\end{equation*}

Recall that $p\geq 5$. By \cite[P.302, Table 4.1]{Gore-82}, $\Mult(M_{11})=\Mult(M_{23})=1$, $\Mult(\A_7)=\ZZ_6$, $\Mult(\A_p)=\ZZ_2$ for $p\not=7$, and $\Mult(\PSL(2,11))=\ZZ_2$. Furthermore,
$\Mult(\PSL(d,q))=\ZZ_{(d,q-1)}$ for $(d,q)\neq(2,4),(2,9),(3,2)$ or $(3,4)$, and for $(d,q)=(2,4),(2,9),(3,2)$ or $(3,4)$, we have $(d,q)=(2,4)$ or $(3,2)$ as $(q^d-1)/(q-1)=p$ and $\Mult(\PSL(2,4))=\Mult(\PSL(3,2))=\ZZ_2$. In all cases, $p\nmid |\Mult(T_2)|$.

Since $N\cong\ZZ_p$ and $X=N.T_2$, by Proposition~\ref{mult}, $X$ is a central extension of $N$ by $T_2$, $X=NX'$ and $X'=H.T_2$, where $H\leq\Mult(T_2)\cap N$. Since $p\nmid |\Mult(T_2)|$, $\Mult(T_2)\cap N=1$ and hence $H=1$. It follows that $H'=T_2$ and $X=N\times H'$. Thus, we may write $X=N\times T_2$.

Let $Y/N=T_1$ with $N\leq Y\unlhd M$. Since $T_1$ is faithful on $\Delta_1$ and $|\Delta_1|=p$, a similar argument to $T_2$ above gives rise to $Y=N\times T_1$. It follows that $M=N\times T_1\times T_2$.

Clearly, $T_1\times T_2$ is characteristic in $M$, an hence $T_1\times T_2\unlhd A$ as $M\unlhd A$. Since $p^3\nmid |T_1\times T_2|$, $T_1\times T_2$ is intransitive on $V(\Ga)$, and by Proposition~\ref{redu}, $T_1\times T_2$ is semiregular, which is impossible. Thus, $t=1$, as required.

\medskip
Now we may let $M/N=T$ for a non-abelian simple group $T$. Since $M/N$ is transitive on $V(\Ga_N)$, $T$ has a subgroup of index $p^2$. Again by Proposition~\ref{simple-ind}, we have
$$T\cong \A_{p^2}, \mbox{ or } \PSL(d,q) \mbox{ with } p^2=(q^d-1)/(q-1) \mbox{ and } d \mbox{ a prime}.$$ For the latter, it is easy to see that $(d,q-1)\not=p$.  By \cite[P.302, Table 4.1]{Gore-82}, we have $p\nmid |\Mult(T)|$, and by a similar argument to $T_2$ above, we may write $M=N\times T$.

Clearly, $T$ is characteristic in $M$ and hence $T\unlhd A$. Since $|V(\Ga)|=p^3$, $T$ cannot be semiregualr on $V(\Ga)$, and by Proposition~\ref{redu}, $T$ is transitive on $V(\Ga)$. Thus, $p^3\mid |T|$.

Assume $T\cong \PSL(d,q)$ with $p^2=(q^d-1)/(q-1)$ and $d$ a prime. Clearly, $p\nmid q$.
If $p\mid (q^d+1)$ then $p$ is divisor of $(q^d-1,q^d+1)$, forcing $p=2$, a contradiction. Then $p\nmid (q^d+1)$. In particular, $p^3\nmid |\PSL(2,q)|$, and since $p^3\mid |T|$, $d$ is an odd prime. Since $p^2=(q^d-1)/(q-1)=q^{d-1}+\cdots+q+1>q^{d-1}$, we have $p>q^{\frac{d-1}{2}}$, forcing $p\nmid (q^i-1)$ for every $1\leq i\leq \frac{d-1}{2}$.

Let $p\mid (q^j-1)$ for some $1+\frac{d-1}{2}\leq j\leq d-1$. Since $p\mid (q^d-1)$, $p$ is a divisor of $(q^d-1)-(q^j-1)=q^j(q^{d-j}-1)$, and since $p\nmid q^j$, $p$ is a divisor of $(q^{d-j}-1)$, which is impossible because $d-j\leq d-(1+\frac{d-1}{2})=\frac{d-1}{2}$.

It follows that $p\nmid (q^i-1)$ for every $1\leq i\leq d-1$. Since $p^2=(q^d-1)/(q-1)$, we have $p^3\nmid |\PGL(d,q)|$, and in particular, $p^3\nmid |\PSL(d,q)|$, contradicting $p^3\mid |T|$.

Thus, $T\cong \A_{p^2}$. Furthermore, $|T:T_\a|=p^2$ for $\a\in V(\Ga_N)$, and $T_\a\cong \A_{p^2-1}$, where $\a$ is an orbit of $N$ on $V(\Ga)$ and $|\a|=p$. Then $T_\a$ has a subgroup of index $p$, which is impossible. This completes the proof of Claim~1.

\medskip
By Claim~1, $M/N$ is abelian, and so $M/N$ is regular on $V(\Ga_N)$. It follows that $M/N\cong\ZZ_p^2$ and $|M|=p^3$. Since $M$ is transitive on $V(\Ga)$, $M$ is regular on $V(\Ga)$. Then we may let $\Ga=\Cay(M,S)$ be a normal Cayley graph, where $R(M)$ is identified with $M$. By Proposition~\ref{Aut-Cay}, $A_1=\Aut(M,S)$, and by Proposition~\ref{nor-2-geo}, $S=\cup_{s\in S}\la s\ra^*$, implying $(p -1)\mid |S|$. Since $N\unlhd A$, $A_1$ fixes $N$ setwise.
Set $\ell=|S|/(p-1)$. Since $\Ga$ is a cover of $\Ga_N$, we have $\Ga_N=\Cay(M/N,\overline{S})$, where $\overline{S}=\{sN\ |\ s\in S\}$, and $M/N\unlhd A/N$. Thus, the stabilizer $(A/N)_N\leq \Aut(M/N,\overline{S})\leq \Aut(M/N)\cong \GL(2,p)$, and noting that $A_1\cong (A/N)_N$, we may write  $$A_1\leq \GL(2,p),\ \mbox{ and } A/N\leq N_{\Aut(\Ga_N)}(M/N),$$
where $N_{\Aut(\Ga_N)}(M/N)$ is the normalizer of $M/N$ in $\Aut(\Ga_N)$. Since $\Ga$ is connected, $\langle S\rangle=M$, and by Proposition~\ref{nor-2-geo}, $M$ is generated by elements of order $p$. By Proposition~\ref{group-p^3}, $M=\ZZ_p^3$ or $\langle a,b,c\ |\ a^p=b^p=c^p=[a,c]=[b,c]=1,c=[a,b]\rangle$.

\medskip
\noindent{\bf Claim~2:} Assume that for all $x,y\in S$ with $\la x\ra\not=\la y\ra$, $x$ is not adjacent to some vertex in $\la y\ra^*$. Then  $\ell(\ell-1)(p-1)n\mid |A_1|$, where $n\geq 1$ and  $\ell=|S|/(p-1)$.

Denote by $\ov{y}_{x+}$ the set of vertices of $\la y\ra^*$ that are adjacent to $x$ in $\Ga$, and by $\ov{y}_{x-}$ the set of vertices of $\la y\ra^*$ that are not adjacent to $x$ in $\Ga$.
Then $|\ov{y}_{x+}|+|\ov{y}_{x-}|=p-1$. Let $\ov{S}=\{\la u\ra^*\ |\ u\in S\}$. Then $|\ov{S}|=\ell$. Let $z\in S$ with $\la z\ra\not=\la x\ra$. By hypothesis, we may take $y^i\in \ov{y}_{x-}$ and $z^j\in\ov{z}_{x-}$. Then $(x,1,y^i)$ and $(x,1,z^j)$ are two
2-geodesics,
and since $\Ga$ is $2$-geodesic transitive, $A_1$ is $2$-transitive on $\ov{S}$. Thus, we may assume that
$$|\ov{y}_{x-}|=n>0,\ \mbox{ for all } x,y\in S \mbox{ with } x\not\in\la y\ra.$$
Clearly, $\Ga(x)\cap \Ga(1)=\cup_{y\in S, y\not\in\la x\ra}\ov{y}_{x+}$, and $|S|=1+|\Ga(x)\cap \Ga(1)|+|\Ga(x)\cap \Ga_2(1)|$ as $\Ga$ has valency $|S|$. Thus, $|\Ga(x)\cap \Ga(1)|=(|\ov{S}|-1)|\ov{y}_{x+}|+p-2=(\ell-1)(p-1-n)+p-2$ and $(p-1)\ell=1+(\ell-1)(p-1-n)+p-2+|\Ga(x)\cap \Ga_2(1)|$, that is, $|\Ga_2(1)\cap \Ga(x)|=(\ell-1)n$. Since $\Ga$ is $2$-geodesic transitive, $A_1$ is transitive on $S$ and $A_{1x}$ is transitive on $\Ga(x)\cap \Ga_2(1)$, which imply that $\ell(\ell-1)(p-1)n\mid |A_1|$, as claimed.

\medskip

Recall that $M=\ZZ_p^3$ or $\langle a,b,c\ |\ a^p=b^p=c^p=[a,c]=[b,c]=1,c=[a,b]\rangle$.

\medskip
\noindent{\bf Case 1:} $M=\ZZ_p^3$.

Let $\Aut(\ZZ_p^3)=\GL(3,p)$ and $Z=Z(\GL(3,p))$. Since  $A_1=\Aut(M,S)\leq \GL(3,p)$ and $S=\cup_{s\in S}\la s\ra^*$, we have $Z\leq A_1\leq \GL(3,p)$. Let $x,y\in S$ with $\la x\ra\not=\la y\ra$. If $x$ is adjacent to all vertices of $\la y\ra$, then every vertex in $\la x\ra$ is adjacent to every vertex in  $\la y\ra$ because $Z\leq A_1$. By Lemma~\ref{clique}, $[\la x\ra\la y\ra]$ is a complete graph, and hence $\la x,y\ra^*\subseteq S$. Since $\Ga$ and $\Ga_N$ have the same valency, Eq.~(\ref{Eq3}) implies that $\Ga$ has valency at most $p^2-1$, and hence $S=\la x,y\ra^*$, contradicting $\la S\ra=M$. Thus, $x$ is not adjacent to some vertex in $\la y\ra^*$. By Claim~2, $\ell(\ell-1)(p-1)n\mid |A_1|$.

Since $M=\ZZ_p^3=\la S\ra$, we have $\ell=|S|/(p-1)\geq 3$. If $\ell=3$,
by Example~\ref{Hamming} we have $\Ga=H(3,p)$, and since $\Ga$ and $\Ga_N$ have the same valency, $\Ga$ has valency $3(p-1)$, which is impossible by Eq.~(\ref{Eq3}). Thus, we may assume that $\ell\geq 4$, and again by Eq.~(\ref{Eq3}), $\Ga_N\cong\K_{p^2}$ or $\ov{H(2,p)}$, and $\Ga$ has valency $p^2-1$ or $(p-1)^2$, respectively.

Suppose $\Ga_N\cong \ov{H(2,p)}$. By Example~\ref{Hamming},  $\Aut(\Ga_N)\cong \Aut(H(2,p))\cong (\SyG_p\times \SyG_p)\rtimes \ZZ_2$, and since $M/N\cong \ZZ_p^2$, we have $A/N\leq N_{\Aut(\Ga_N)}(M/N)\lesssim (\AGL(1,p)\times \AGL(1,p))\rtimes \ZZ_2$. It follows that
$|(A/N)_N|\mid  2(p-1)^2$ and hence $|A_1|\mid 2(p-1)^2$.
Since $\ov{H(2,p)}$ has valency $(p-1)^2$, we have $|S|=(p-1)^2$ and $\ell=p-1$. Then $(p-1)^2(p-2)n\mid |A_1|$ and so $(p-1)^2(p-2)n\mid 2(p-1)^2$, a contradiction.

Suppose $\Ga_N\cong \K_{p^2}$. Then $|S|=p^2-1$ and $\ell=p+1$. Since $p(p^2-1) \mid |A_1|$ and $A_1\leq \GL(2,p)$, we have $|\GL(2,p):A_1|\leq p-1$, and by Lemma~\ref{autEp^3}~(1), $A_1$ contains $\GL(2,p)'\cong\SL(2,p)$.
Since $A_1=\Aut(M,S)$ fixes $N$ setwise, $A_1$ is a subgroup of $\Aut(M)_N$, the block stabilizer of $N$ in $\Aut(M)$. We may assume that $M=\ZZ_p^3=\la a\ra\times\la b\ra\times\la c\ra$ and $N=\la c\ra$. By Lemma~\ref{EAp3C}~(3), $\Aut(M)_N$ has one conjugate class of subgroup isomorphic to $\SL(2,p)$, of which one has orbit set $\{\{c^i\},c^i\la a,b\ra^*\mid i\in\ZZ_p\}$.
Thus, each orbit of $A_1$ on $V(\Ga)$ is union of some orbits of $\SL(2,p)$.
Since $A_1$ is transitive on $S$, $S$ is an orbit of $A_1$ with length $p^2-1$, and hence $S=c^i\la a,b\ra^*$ for some $i\in\ZZ_p$.
Since $\la S\ra=M$, $c^i\neq1$ and $\la ch \ra^*\subseteq S$ for all $h\in\la a,b\ra^*$. Thus, $|S|>p^2-1$, a contradiction.

\medskip
\noindent{\bf Case 2:} $M=\langle a,b,c\ |\ a^p=b^p=c^p=[a,c]=[b,c]=1,c=[a,b]\rangle$.

In this case, $N=\langle c\rangle\unlhd A$. By Eq.~(\ref{Eq3}),
$\Ga$ has valency at least $2(p-1)$. Then $|S|\geq 2(p-1)$, and since $A_1=\Aut(M,S)$ fixes $\la c\ra$ setwise, the arc-transitivity of $A$ implies $\la c\ra\cap S=\emptyset$. Thus, $\langle x,y\rangle=M$ and $[x,y]\not=1$, for all $x,y\in S$ with $x\not\in \la y\ra$.

Suppose that $x$ is adjacent to all elements in $\la y\ra^*$. Then $xy^i\in S$ for all $i\in\ZZ_p^*$.
If $x^j\in \la x\ra^*$ is not adjacent to $y^k\in\la y\ra^*$, then $A_1$ has an element $\a$ interchanging $x^j$ and $y^k$ and hence  $\la x\ra^*$ and $\la y\ra^*$ because $(x^j,1,y^k)$ and $(y^k,1,x^j)$ are $2$-geodesics, implying that $x^\a\in\la y\ra^*$ is adjacent to all vertices in $\la x\ra^*$. Thus, we may always assume that $\la x\ra^*$
has an element $x^s$ adjacent to some $y^t\in\la y\ra^*$ with $1\not=s\in\ZZ_p^*$ and $t\in\ZZ_p^*$. Thus, $x^sy^{-t}\in S$ and so $\la x^sy^{-t}\ra^*\in S$. Then
$(x^sy^{-t})^{s^{-1}}=xy^{-s^{-1}t}[y^{-t},x^s]^{2^{-1}s^{-1}(s^{-1}-1)}
=xy^{-s^{-1}t}c_1\in S$, where $c_1=[y^{-t},x^s]^{2^{-1}s^{-1}(s^{-1}-1)}\in \la c\ra$. Since $s\not=1$ and $[x,y]\not=1$, we have $c_1\not=1$. Recall that $xy^{-s^{-1}t}\in S$. Then we have two edges $\{1,xy^{-s^{-1}t}\}$ and $\{1,xy^{-s^{-1}t}c_1\}$, implying the subgraph $[N, xy^{-s^{-1}t}N]$ of $\Ga$ is not a matching, contradicting that $\Ga$ is a cover of $\Ga_N$. Thus, $x$ is not adjacent to some vertex in $\la y\ra^*$. By Claim~2,  $|\Ga_2(1)\cap \Ga(x)|=(\ell-1)n$ and $\ell(\ell-1)(p-1)n\mid |A_1|$, and by Eq.~(\ref{Eq3}), $\Ga_N\cong H(2,p)$, $\K_{p^2}$ or $\ov{H(2,p)}$.

Let $\Ga_N\cong H(2,p)$. Then $\Ga$ has valency $2(p-1)$, and since $S=\cup_{s\in S}\la s\ra^*$, we may let $S=\la x\ra^* \cup \la y\ra^*$. Thus, $M=\langle S\rangle=\la x,y\ra$, and so $M$ has an automorphism mapping $x$ to $a$ and $y$ to $b$. It follows that $\Ga\cong \Cay(M,\la a\ra^*\cup \la b\ra^*)$, which corresponds to the first graph of Theorem~\ref{Thm-p^3}~(2)-(iv).

Let $\Ga_N\cong \ov{H(2,p)}$. Then $\ell=(p-1)^2/(p-1)=p-1$, and a similar argument to Case 1 yields that $|A_1|\mid 2(p-1)^2$. Since $(p-1)^2(p-2)n\mid |A_1|$, we have $(p-1)^2(p-2)n\mid 2(p-1)^2$, a contradiction.

Let $\Ga_N\cong \K_{p^2}$. Then $|S|=p^2-1$ and $\ell=p+1$. Since $\ell(\ell-1)(p-1)n\mid |A_1|$, we have $p(p^2-1) \mid |A_1|$ and $|\GL(2,p):A_1|\leq p-1$ as $A_1\leq \GL(2,p)$. By Lemma~\ref{autEp^3}, $\Aut(M)$ has one conjugate class of subgroup isomorphic to $A_1$, and $\GL(2,p)'=\SL(2,p)\leq A_1$. By Theorem~\ref{NPp3-2}, we may let $$A_1\leq \Aut(M,R)\cong\GL(2,p), \mbox{ where } R=\la b\ra^*\cup_{i\in\ZZ_p}\la b^iab^i\ra^* \mbox{ with } |R|=p^2-1.$$
Let $|\Aut(M,R):A_1|=m<p-1$. By Theorem~\ref{NPp3-2}~(1), $A_1$ has only one orbit with length $p^2-1$, that is $R$, and so $S=R$, which corresponds the second graph of  Theorem~\ref{Thm-p^3}~(2)-(iv).
Now let $|\Aut(M,R):A_1|=p-1$. Then $|A_1|=p(p^2-1)$.
By Theorem~\ref{NPp3-2}~(1), $A_1$ has $p$ orbits of length $1$, and $p$ orbits of length $p^2-1$, witch are $c^iR$, $i\in\ZZ_p$, by Theorem~\ref{NPp3-2}~(2). It follows that $S=c^iR$ for some $i\in \ZZ_p$. Note that $SS=S^{-1}S\subseteq \{1\}\cup \Ga(1)\cup \Ga_2(1)$. Suppose $c^i\not=1$. By Theorem~\ref{NPp3-2}~(3), $|S^{-1}S|=|(c^iR)^{-1}c^iR|=|RR|=p^2+(p^2-1)(p-1)$, and then $|\Ga_2(1)|\geq p^2+(p^2-1)(p-1)-p^2=(p^2-1)(p-1)$,
forcing $|\Ga_2(1)|=p(p^2-1)$
because $|A_1|=p(p^2-1)$ and $A_1$ is transitive on $\Ga_2(1)$. However, $|V(\Ga)|\geq 1+|\Ga(1)|+|\Ga_2(1)|=p^3+p^2-p>p^3$, a contradiction.
Thus, $c^i=1$ and $S=R$, which corresponds Theorem~\ref{Thm-p^3}~(2)-(iv).\qed

By Theorem~\ref{Thm-p^3}, the graph defined in Theorem~\ref{NPp3-2} is not isomorphic to any other graphs in Theorem~\ref{Thm-p^3}. Then the proof of Theorem~\ref{Thm-p^3} implies the following.

\begin{corollary}\label{coro1}
Let $\Ga$ be the graph defined in Theorem~\ref{NPp3-2}.
Then $\Ga$ is a normal Cayley graph and $\Aut(\Ga)\cong R(E(p^3))\rtimes \GL(2,p)$.
\end{corollary}

\end{document}